\documentclass[12pt,reqno]{amsart}
\usepackage{amsmath,amsfonts,amssymb,amscd,amsthm,amsbsy,epsf}

\usepackage{ucs}
\usepackage[utf8]{inputenc}
\usepackage{mathabx}
\usepackage{caption}
\usepackage[UKenglish]{babel}
\usepackage[dvips]{graphicx}
\graphicspath{{}}

\usepackage{hyperref}
\usepackage{comment}
\usepackage{multirow}
\usepackage{color}

\usepackage{longtable}
\allowdisplaybreaks

\title[KAM quasi-periodic tori for the dissipative spin-orbit problem]
{KAM quasi-periodic tori for the dissipative spin-orbit problem}

\author[R. Calleja]{Renato Calleja}

\address{ Department of Mathematics and Mechanics, IIMAS, National
  Autonomous University of Mexico (UNAM), Apdo. Postal 20-126,
  C.P. 0100, Mexico D.F. (Mexico)}

\email{celleja@mym.iimas.unam.mx}

\author[A. Celletti]{Alessandra Celletti}

\address{ Department of Mathematics, University of Rome Tor Vergata,
  Via della Ricerca Scientifica 1, 00133 Rome (Italy)}

\email{celletti@mat.uniroma2.it}

\author[J. Gimeno]{Joan Gimeno}

\address{ Department of Mathematics, University of Rome Tor Vergata,
  Via della Ricerca Scientifica 1, 00133 Rome (Italy)}

\email{gimeno@mat.uniroma2.it}

\author[R. de la Llave]{Rafael de la Llave}

\address{ School of Mathematics, Georgia Institute of Technology, 686
  Cherry St.. Atlanta GA. 30332-0160 (USA) }

\email{rafael.delallave@math.gatech.edu}

\thanks{R.C. was partially supported by UNAM-DGAPA PAPIIT Project IN
  101020.  A.C. was partially supported the MIUR Excellence Department
  Project awarded to the Department of Mathematics, University of Rome
  Tor Vergata, CUP E83C18000100006, EU H2020 MSCA ETN
  Stardust-Reloaded Grant Agreement 813644, MIUR-PRIN 20178CJA2B ``New
  Frontiers of Celestial Mechanics: theory and
  Applications''. J.G. has been supported by the Spanish grants
  PGC2018-100699-B-I00 (MCIU/AEI/FEDER, UE), the Catalan grant 2017
  SGR 1374 and MIUR-PRIN 20178CJA2B ``New Frontiers of Celestial
  Mechanics: theory and Applications''.  J.G. thanks the School of
  Mathematics of GT for its hospitality in Spring 2019 and Fall
  2019. R.L. has been partially supported by NSF grant DMS 1800241.}

\date{\today}

\newcommand{\re}[1]{{#1}_\mathtt{R}}
\newcommand{\I}{\mathtt{i}}
\newcommand{\im}[1]{{#1}_\mathtt{I}}
\let\Re\relax \DeclareMathOperator{\Re}{Re}
\let\Im\relax \DeclareMathOperator{\Im}{Im}

\newcommand{\ve}[1]{\boldsymbol{#1}}
\newcommand{\eps}{\varepsilon}
\newcommand{\dis}{\eta}
\newcommand{\ecc}{e}
\newcommand{\omg}{\omega}
\newcommand{\param}{e}

\newtheorem{thm}{Theorem}[section]
\newtheorem{meta-thm}[thm]{Meta-Theorem}

\newtheorem{rem}[thm]{Remark}

\newtheorem{defn}[thm]{Definition}

\newtheorem{lem}[thm]{Lemma}
\newtheorem{alg}[thm]{Algorithm}

\newcommand\beq[1]{ \begin{equation}\label{#1} }
\newcommand{\eeq}{ \end{equation} }

\newcommand\beqa[1]{ \begin{eqnarray} \label{#1}}
\newcommand{\eeqa}{ \end{eqnarray} }
\newcommand{\beqano}{ \begin{eqnarray*} }
\newcommand{\eeqano}{ \end{eqnarray*} }
\newcommand\equ[1]{{\rm (\ref{#1})}}
\newcommand{\R}{\mathbb{R}}

\newcommand{\red}{\textcolor{red}}

\renewcommand{\thesection}{\S\arabic{section}}
\renewcommand{\thethm}{\arabic{section}.\arabic{thm}}

\begin{document}
\begin{abstract}
We provide evidence of the existence of KAM quasi-periodic attractors
for a dissipative model in Celestial Mechanics. We compute the
attractors extremely close to the breakdown threshold.

We consider the spin-orbit problem describing the motion of a triaxial
satellite around a central planet under the simplifying assumption
that the center of mass of the satellite moves on a Keplerian orbit,
the spin-axis is perpendicular to the orbit plane and coincides with
the shortest physical axis. We also assume that the satellite is
non-rigid; as a consequence, the problem is affected by a dissipative
tidal torque that can be modeled as a time-dependent friction, which
depends linearly upon the velocity.

Our goal is to fix a frequency and compute the embedding of a smooth
attractor with this frequency. This task requires to adjust a drift
parameter.

We have shown in \cite{CCGL20a} that it is numerically efficient to
study Poincar\'e maps; the resulting \emph{spin-orbit map} is
conformally symplectic, namely it transforms the symplectic form into
a multiple of itself. In \cite{CCGL20a}, we have developed an
extremely efficient (quadratically convergent, low storage
requirements and low operation count per step) algorithm to construct
quasi-periodic solutions and we have implemented it in extended
precision. Furthermore, in \cite{CCL20} we have provided an
``a-posteriori'' KAM theorem that shows that if we have an embedding
and a drift parameter that satisfy the invariance equation up to an
error which is small enough with respect to some explicit condition
numbers, then there is a true solution of the invariance
equation. This a-posteriori result is based on a Nash-Moser hard
implicit function theorem, since the Newton method incurs losses of
derivatives.

The goal of this paper is to provide numerical calculations of the
condition numbers and verify that, when they are applied to the
numerical solutions, they will lead to the existence of the torus for
values of the parameters extremely close to the parameters of
breakdown.  Computing reliably close to the breakdown allows to
discover several interesting phenomena, which we will report in
\cite{CCGL20c}.

The numerical calculations of the condition numbers presented here are
not completely rigorous, since we do not use interval arithmetic to
estimate the round off error and we do not estimate rigorously the
truncation error, but we implement the usual standards in numerical
analysis (using extended precision, checking that the results are not
affected by the level of precision, truncation, etc.). Hence, we do
not claim a computer-assisted proof, but the verification is more
convincing that standard numerics. We hope that our work could
stimulate a computer-assisted proof.
\end{abstract}


\maketitle


\section{Introduction}
Kolmogorov-Arnold-Moser (hereafter KAM) theory
(\cite{Kolmogorov54,Arnold63a,Moser62}) concerns the existence of
quasi-periodic motions in non-integrable dynamical systems.  In its
original formulation, it was applied to nearly-integrable Hamiltonian
systems.

An important recent development is the \emph{a-posteriori KAM theory}
(see \cite{LlaveGJV05,Llave01c}) that does not require that the system
is close to integrable, but rather that there is an approximate
solution of an invarance equation that satisfies some non-degeneracy
conditions.  Given an a-posteriori KAM theorem, one does not need to
justify the way that the approximate solution is constructed (it could
be done by formal expansions or just by numerical tries), but one must
provide rigorous estimates on the error of the invariance equation and
the condition numbers involved in the theorem statement.

The KAM theory has been extended to general systems (see, e.g.,
\cite{Moser67}). This theory fixes the frequency of the quasi-periodic
orbit searched, but adjusting parameters in the system. This general
KAM theory is even more effective if the system preserves some
geometric structures (\cite{BroerHTB90,BroerHS96,broer2005}). From the
mathematical point of view, the number of parameters to adjust may be
reduced (e.g., in the Hamiltonian case, there are no parameters to be
adjusted). Numerically, one can use identities coming from the
geometry to develop fast algorithms that also require small storage
space and enjoy good stability properties.  For the purposes of our
paper, the most relevant development is \cite{CallejaCL11}, which
established an a-posteriori KAM theorem and presented efficient
numerical algorithms for \emph{conformally symplectic systems} (that
is, systems that transform the symplectic form into a multiple of
itself). Conformally symplectic systems appear in a variety of
applications, including Euler-Lagrange equations of exponentially
discounted Lagrangians, thermostats, etc.

The goal of this paper is to study the applicability of a-posteriori
KAM theory for a specific model of Celestial Mechanics known as the
\emph{spin-orbit problem with tidal torque}. This model describes the
rotational motion of a non-rigid triaxial ellipsoid orbiting around a
point-mass planet. We assume that the planet moves in a Keplerian
orbit, the rotation axis is perpendicular to the orbital plane and
aligned with the shortest physical axis of the satellite. Furthermore,
we assume that the system experiences a tidal force proportional to
the velocity, which makes it into a conformally symplectic
system. This model has been studied in
\cite{CellettiC2009,Massetti19,Locatelli}.

Efficient numerical methods to find quasi-periodic orbits in the
spin-orbit model were implemented in \cite{CCGL20a}. Taking advantage
of the extreme efficiency of the methods, modern programming tools and
the power of modern hardware. The calculations of \cite{CCGL20a} were
run in high precision and produced the parameterization of
quasi-periodic orbits and adjusted parameters that solve the
invariance equations with very high accuracy, even very close to the
breakdown\footnote{As a matter of fact, there is no alternative
  numerical method that can compute as close to the breakdown, so that
  the estimates of this paper are the best estimates for the
  threshold, since the solutions we can compute have all the signs of
  being very deteriorated.}.

The goal of this paper is to study the application of the a-posteriori
theorem in \cite{CallejaCL11} to the calculations in
\cite{CCGL20a}. We take the calculations in \cite{CCGL20a}, and
evaluate numerically the condition numbers required in
\cite{CallejaCL11}. Similar results for an explicitly given mapping
appear in \cite{CCL20}. In the present problem, the map considered is
not given by an explicit formula, but is obtained by integrating an
ordinary differential equation.  This requires new analysis and
numerical studies of the variational equations.

The results presented here come short of a full computer-assisted
proof, since the evaluation of the error and the condition numbers are
not completely rigorous. We do not take into account round-off or
truncation errors.

We certainly hope that the present effort could serve as inspiration
for others to close the gap and provide a true computer-assisted proof
and, needless to say, we would be happy to provide detailed data and
encouragement.  Even if not the final word on existence, we think that
the work presented goes beyond the regular standards of numerical
computations and is a significant progress in the area of the
computations of tori, even close to the breakdown. We think that it is
rather remarkable that the algorithms inspired by the theory are also
the most efficient ones.

Computing close to the breakdown and being able to trust the
computation is not just an affectation, but uncovers new phenomena
that present a challenge to mathematics.

We note that, even if the computation is doable, but delicate for
values of the perturbation close to the threshold, it remains
extremely reliable and easy for many values of astronomical interest,
so that KAM theory and their algorithms become a relevant tool to
astronomers, overcoming the concerns --relevant at the time they were
written-- of \cite{Henon66}.


\vskip.1in

This paper is organized as follows. The equation of motion describing
the dissipative spin-orbit problem is shortly recalled in
Section~\ref{sec:model}. We study the Poincar\'e map associated to
such a model in Section~\ref{sec:map}; in this way we obtain a {\sl
  spin-orbit map}, which is conformally symplectic, and we compute the
corresponding conformally symplectic factor, which is the term by
which the symplectic form gets multiplied, when the map is applied to
the the symplectic form. Then, we use the KAM theorem for conformally
symplectic maps formulated in \cite{CCL20} (see
Section~\ref{sec:KAM}). Contrary to the implementation to the standard
map, the application of the theorem to the spin-orbit problem is more
complex and it requires a careful computation of some constants as
described in Section~\ref{sec:Qestimates}. This procedure leads to the
final results that we present in Section~\ref{sec:estimates} for two
different frequencies: the golden ratio and a second frequency between
one and the golden ratio.

\section{The spin-orbit problem with tidal torque}
\label{sec:model}
For the sake of motivation, in this section we present the physical
basis of the model considered. Even if this motivates the questions
asked, it is logically independent of the analysis.

Consider the motion of a non-rigid satellite $\mathcal{S}$ that we
assume to have a triaxial shape and principal moments of inertia
$\mathcal{A} < \mathcal{B} < \mathcal{C}$. We assume that the
barycenter of the satellite $\mathcal{S}$ moves on an elliptic
Keplerian orbit with semimajor axis $a$, eccentricity $\ecc$, and with
the planet $\mathcal{P}$ in one focus.  The satellite rotates around
the smallest physical axis, in such a way that the spin-axis is
perpendicular to the orbit plane (see, e.g.,
\cite{Beletsky,Celletti90I,Celletti2010,LaskarC,Wisdom}).

We normalize the units of measure of time so that the orbital period
$T_{orb}$ is equal to $2\pi$, which implies that the mean motion is
$n=2\pi/T_{orb}=1$; we introduce the {\sl perturbative parameter}
$\varepsilon$, which measures the equatorial ellipticity of the
satellite:
\begin{equation}
 \label{eq.eps}
 \eps \coloneq \frac{3}{2} \frac{\mathcal{B}-\mathcal{A}}{\mathcal{C}}\ .
\end{equation}
We denote by $x$ the angle between the largest physical axis of the
triaxial satellite and the periapsis line. The equation of motion of
the spin-orbit problem, using the formulation in
\cite{Macdonald,peale} for the tidal torque, is given by
\begin{equation}\label{eq.spin-diss}
 \frac{d^2x(t)}{dt^2} + \eps \biggl(\frac{a}{r(t)}\biggr)^3 \sin \bigl(2 x(t) -
 2f(t)\bigr) =  - \dis
 \biggl(\frac{a}{r(t)}\biggr)^6 \biggl(\frac{dx(t)}{dt} - \frac{d f(t)}{dt}
 \biggr)\ ,
\end{equation}
where $r(t)=r(t;\ecc)$ and $f(t)=f(t;\ecc)$ are the orbital radius and
the true anomaly of the Keplerian ellipse, and $\dis > 0$ is the {\sl
  dissipative constant} depending on the physical features of the
satellite. Denoting by $u$ the eccentric anomaly, then
\[
r=a(1-\ecc \cos u)\ ,\qquad
\tan \biggl(\frac{f}{2}\biggr) = \sqrt{{{1+\ecc}\over{1-\ecc}}} \tan
\biggl(\frac{u}{2} \biggr)\ .
\]
For $\dis=0$ the model becomes conservative and takes a
nearly-integrable form with $\eps$ being the {\sl perturbing
  parameter}.  We also introduce the spin-orbit problem with tidal
torque \sl averaged \rm over one orbital period (see, e.g.,
\cite{peale, CCGL20a}) as given by the equation
\begin{equation}
 \label{eq.avg-spinxy}
\frac{d^2x(t)}{dt^2} + \eps \Big({a\over r(t)}\Big)^3 \sin
\bigl(2x(t)-2f(t) \bigr) = -\dis \bar
L(\ecc)\biggl(\frac{dx(t)}{dt}-{{\bar
    N(\ecc)}\over {\bar L(\ecc)}}\biggr)\ ,
\end{equation}
where
\begin{equation*}
 \begin{split}
\bar L(\ecc)&\coloneq {1\over{(1-\ecc^2)^{9/2}}} \biggl(1+3\ecc^2+{3\over
  8}\ecc^4 \biggr) \ , \\
\bar N(\ecc)&\coloneq
     {1\over{(1-\ecc^2)^6}} \biggl(1+{{15}\over 2}\ecc^2+{{45}\over
       8}\ecc^4+ {5\over {16}}\ecc^6 \biggr)\ .
 \end{split}
\end{equation*}

\section{The conformally symplectic spin-orbit map}
\label{sec:map}

Following \cite{CCGL20a}, we introduce a discrete system, which is
obtained by computing the Poincar\'e map $P_e$ associated to
\equ{eq.spin-diss}. Precisely, setting $y=\dot x$, we can write the
map as
\begin{equation}
\label{eq.Pe}
P_e(x_0,y_0;\eps) \coloneq
\begin{pmatrix}
 x(2\pi;x_0,y_0,\eps) \\ y(2\pi;x_0,y_0,\eps)
\end{pmatrix}\ ,
\end{equation}
where $x(2\pi;x_0,y_0,\eps) $ and $y(2\pi;x_0,y_0,\eps)$ denote the
solution of \equ{eq.spin-diss} at time $t=2\pi$ with initial
conditions $(x_0,y_0)$ at $t=0$. Writing $P_e$ in components, say $P_e
\equiv (P_e^{(1)},P_e^{(2)})$, the {\sl spin-orbit Poincar\'e map}
becomes:
\begin{equation}
\label{SOmap}
 \begin{split}
\bar x&= P_e^{(1)}(x,y;\eps)\ ,\\
\bar y&= P_e^{(2)}(x,y;\eps)\ .
 \end{split}
\end{equation}
For numerical reasons, it is better to consider the change of
coordinates
\begin{equation}
 \label{eq.xy2bg}
 \Psi _\ecc \coloneq 2\pi
 \begin{pmatrix}
  1 & 0 \\ 0 & 1 - \ecc
 \end{pmatrix}
\end{equation}
and define the map $G _\ecc \coloneq \Psi _\ecc \circ P _\ecc \circ
\Psi _\ecc^{-1}$ which can be computed accurately by numerical
integrators such as \cite{HairerNW1993,JorbaZ2005}.

The map \equ{SOmap}, equivalently $G_\ecc$, inherits several
properties of the continuous system \equ{eq.spin-diss}. In particular,
the map is {\sl conformally symplectic}, which means that it
transforms the symplectic form into a multiple of itself, according to
the following definition.

\begin{defn}\label{def:conformallysymplectic}
Let $\mathcal{M} = \mathbb{T}^n\times U$ with $U\subseteq
\mathbb{R}^n$ an open and simply connected domain with smooth
boundary. We endow $\mathcal{M}$ with a symplectic form $\Omega$.  A
diffeomorphism $f\colon \mathcal M\rightarrow\mathcal M$ is
conformally symplectic, if there exists a function $\lambda \colon
\mathcal{M}\to\mathbb{R}$ such that
\begin{equation}\label{CS}
f^* \Omega = \lambda\Omega\ ,
\end{equation}
where $f^*$ denotes the pull--back of $f$.
\end{defn}

We will call $\lambda$ the conformal factor. For $\lambda=1$ we have a
symplectic diffeomorphism.  In the following, we will consider the
family $ P _\ecc \colon {\mathcal M}\rightarrow{\mathcal M}$, defined
in \eqref{eq.Pe}, of diffeomorphisms depending on a parameter $\ecc\in
[0,1)$ to which we refer as the {\sl drift parameter}. In this case
  \eqref{CS} is replaced by
\begin{equation}\label{CSmu}
P_ \ecc^* \Omega = \lambda\Omega\ .
\end{equation}

\vskip.1in

The definition of conformally symplectic continuous systems is given
as follows.

\begin{defn}
A vector field $X$ is a conformally symplectic flow if, denoting by
$L_X$ the Lie derivative, there exists a function $\lambda \colon
\mathbb{R}^{2n}\to\mathbb{R}$ such that
\begin{equation}
\label{flow}
 L_X \Omega = \lambda \Omega\ .
\end{equation}
\end{defn}

\vskip.1in

If $\Phi_t$ denotes the flow at time $t$, then \eqref{flow} implies
that
\[
(\Phi_t)^*\Omega=\exp({\lambda t})\Omega\ .
\]

\vskip.1in

The dissipative spin-orbit model \equ{eq.spin-diss} is an example of a
conformally symplectic vector field.  An important result for our
purposes is that the Poincar\'e map associated to a conformally
symplectic vector field is a conformally symplectic map. As a
consequence, the spin-orbit Poincar\'e map defined in \equ{SOmap} is
conformally symplectic with the conformally symplectic factor given by
\begin{equation}
 \label{eq.simplfactor}
 \lambda(x,y) = \sigma|\det D P_ \ecc (x,y;\eps) |, \qquad \sigma = \pm 1\ ,
\end{equation}
where $\sigma$ denotes the orientation of $P _\ecc$.

As shown in \cite{CCGL20a}, the conformal factor is given explicitly
in terms of the orbital eccentricity and the dissipative parameter:
\begin{equation}\label{lambda}
\lambda = \exp \biggr(-\dis \pi \frac{3 \ecc^4+24 \ecc^2+8}{4
\left(1-\ecc^2\right)^{9/2}}\biggl)\ .
\end{equation}
When $\dis>0$ we have a contractive system, if $\dis<0$ we have an
expansive system and if $\dis=0$ we have a symplectic system.  In the
following we will just consider the contractive case with $\dis> 0$.

\section{KAM theorem and invariant attractors} \label{sec:KAM}

The statement of the KAM theorem that we will apply to the spin-orbit
problem requires a set of preliminary notations and notions. We start
to give, in Section~\ref{sec:norms}, the definition of the norms and
some results on Cauchy estimates on the derivatives. In
Section~\ref{sec:cohomology} we give the definition of Diophantine
frequency and we present some results on the solution of the
cohomology equation. The definition of KAM attractor and the
invariance equation to be satisfied is given in
Section~\ref{sec:invariance}. Finally, the statement of the KAM
theorem, borrowed from \cite{CCL20}, is given in
Section~\ref{sec:KAMthm}.

\subsection{Norms and Cauchy estimates}\label{sec:norms}
The norm of a vector ${\underline v}= \left(
\begin{smallmatrix}
 v _1 \\ v _2
\end{smallmatrix}
\right)\in{\mathbb R}^2$ is defined as
\[
\|{\underline v}\| \coloneq |v_1|+|v_2|\ .
\]

The norm of a matrix $A= \left(
\begin{smallmatrix}
  a _{11} & a _{12} \\ a _{21} & a _{22}
\end{smallmatrix}
\right)\in{\mathbb R}^2 \times {\mathbb R}^2$ is defined as
\[
\|A\| \coloneq \max\bigl\{|a _{11}| + |a _{21}|,\ |a _{12}| + |a
_{22}|\bigr\}\ .
\]
Next, we consider the norm of functions and vector functions. To this
end, for $\rho>0$ we introduce the complex extensions of a torus
${\mathbb T}$, a set $B$ and the manifold ${\mathcal M}={\mathbb T}
\times B$ as
\begin{align}
  {\mathbb T}_\rho&\coloneq  \{x+\I y\in{\mathbb C}/{\mathbb
Z}\colon \ x\in{\mathbb  T}\ ,\ |y|\leq\rho\}\ , \label{eq.Trho} \\
  B_\rho&\coloneq \{x+\I y\in{\mathbb C}\colon \ x\in B\ ,\quad |y|\leq \rho\}\ , \nonumber  \\
  {\mathcal M}_\rho&\coloneq {\mathbb  T}_\rho\times B_\rho\ . \nonumber
\end{align}
By ${\mathcal A}_{\rho}$ we denote the set of functions analytic in
the interior of ${\mathbb T}_\rho$ and extending continuously to the
boundary of ${\mathbb T}_\rho$. This set is endowed with the norm
\begin{equation}
\label{eq.normArho} \|f\|_{\rho} \coloneq \sup_{z\in{\mathbb T}_\rho}
|f(z)|\ .
\end{equation}

Similarly, for a vector valued function $f=(f_1,f_2,\dotsc,f_n)$,
$n\geq 1$, we define the norm
\begin{equation}
\label{normv}
 \|f\|_{\rho} \coloneq \|f_1\|_{\rho}+\|f_2\|_{\rho}+ \dotsb
+\|f_n\|_{\rho}\ .
\end{equation}
If $F$ denotes an $n_1\times n_2$ matrix valued function, then we
define its norm as
\begin{equation}
 \label{normm}
 \|F\|_{\rho} \coloneq \sum_{i=1}^{n_1} \sup_{j=1,\dotsc,n_2}
\|F_{ij}\|_{\rho}\ .
\end{equation}

The following classical lemma gives a bound on the derivatives on
smaller domains than the initial function (see, e.g., \cite{CCL20} for
its proof).

\begin{lem} \label{lem:Cauchy}
  Given a function $h\in{\mathcal A}_{\rho}$, its first derivative can
  be bounded as
  \begin{equation} \label{cc}
    \|Dh\|_{\rho-\delta}\leq \delta^{-1}\ \|h\|_\rho\ ,
  \end{equation}
  where $0<\delta<\rho$.
\end{lem}

\subsection{Diophantine frequency and the cohomology equation}
\label{sec:cohomology}

One of the main assumptions in KAM theory is that the frequency
satisfies a Diophantine assumption that, in view of the application of
KAM theory to the spin-orbit map \equ{SOmap}, we introduce as follows.

\begin{defn}
\label{defn.Diophantine}
Let $\omega\in{\mathbb R}$ and let $\tau\geq 1$, $\nu>0$. The number
$\omega$ is said Diophantine of class $\tau$ and constant $\nu$,
$\omega\in{\mathcal D}(\nu,\tau)$, if for all $q \in \mathbb{Z}$ and
$k \in {\mathbb Z}\backslash\{0\}$, it satisfies the following
inequality
\begin{equation}
|\omega \,k-q|\ \geq\ \nu |k|^{-\tau} \ . \label{DC}
\end{equation}
\end{defn}

Another important ingredient at the basis of the proof of the KAM
theorem is the solution of a cohomology equation of the form
\begin{equation} \label{difference}
\varphi(\theta+\omega)-\lambda \varphi(\theta)=\vartheta(\theta)\ ,
\end{equation}
where $\theta\in{\mathbb T}$ and $\vartheta$ is a Lebesgue measurable
function.

The following lemmas yield the existence of a solution of
\equ{difference} given by a Lebesgue measurable function
$\varphi$. The first result, Lemma~\ref{contractive}, is valid when
$|\lambda| \ne 1$ and $\omega\in {\mathbb R}$. It gives an estimate on
the solution which depends on $\lambda$ and indeed explodes as
$|\lambda|$ tends to $1$. The second result, Lemma~\ref{neutral}, is
valid for any $\lambda$ and Diophantine frequency $\omega$. It
provides a uniform estimate of the solution.  We refer to
\cite{CallejaCL11,CCL20} for the proofs of the
Lemmas~\ref{contractive} and \ref{neutral}.

In \cite{CallejaCL11}, one can find also estimates that are uniform
for $\lambda \in [A^{-1}, A]$ for $A> 1$ and, hence allow to study the
(singular) limit of zero dissipation. These estimates are very similar
to the estimates in Lemma~\ref{neutral} (they use the Diophantine
condition and they entail a loss of domain).

\begin{lem}\label{contractive}
Let $|\lambda| \ne1$ and $\omega\in {\mathbb R}$.  Given any Lebesgue
measurable function $\vartheta$, there exists a Lebesgue measurable
function $\varphi$ which satisfies \eqref{difference} and which is
bounded by
\[
\|\varphi\|_{\rho}  \le\big|\,|\lambda|
-1\, \big|^{-1}\|\vartheta\|_{\rho}\ .
\]
The derivatives of $\varphi$ with respect to $\lambda$ are bounded by
\[
\|D_\lambda^j \varphi\|_{\rho} \le {j!\over
  {\bigl|\, |\lambda|-1\, \bigr|^{j+1}}}\ \|\vartheta\|_{\rho}\ ,\qquad j\geq
1\ .
\]
\end{lem}

\begin{lem}\label{neutral}
Assume that $\lambda \in [A_0, A_0^{-1}]$ for some $0 < A_0< 1$ in
\eqref{difference} and let $\omega\in {\mathcal D}(\nu,\tau)$.  Let
$\vartheta \in {\mathcal A}_{\rho}$, $\rho>0$, be a function such that
\[
\int_{{\mathbb T}} \vartheta(\theta)\, d\theta =0\ .
\]
Then, there exists one, and only one, solution of \eqref{difference}
with zero average:
\[
\int_{{\mathbb T}} \varphi(\theta)\, d\theta =0\ .
\]
Moreover, if $\varphi \in {\mathcal A}_{\rho-\delta}$ for $0 <
\delta < \rho$, then we have

\begin{equation} \label{estimate} \|\varphi\|_{\rho-\delta} \le C_0\ \nu^{-1}\
\delta^{-\tau} \|\vartheta\|_{\rho}\ ,
\end{equation}
where
\begin{equation} \label{Cu} C_0={1\over
{(2\pi)^\tau}}\ {\pi\over {2^\tau (1+\lambda)}}\
\sqrt{{\Gamma(2\tau+1)}\over 3}
\end{equation}
and $\Gamma$ denotes the gamma function.
\end{lem}

We remark that \cite{FHL} provides a better estimate for the constant
$C_0$ in the symplectic case. Its expression is more complicated than
\equ{Cu}. However, for our parameter values, it seems that the
estimate \equ{Cu} suffices to reach the final result of getting
analytic estimates close to the break-down.

\subsection{KAM attractor and the invariance equation}
\label{sec:invariance}
In this Section, we introduce the definition of a KAM attractor with
Diophantine frequency $\omega$ for a family $f_\param$ of conformally
symplectic maps. We call $\param$ the drift parameter, since we
recognize that the drift is related to the eccentricity, although the
drift might in principle coincide with a different parameter.  This
will require to satisfy the invariance equation \eqref{invariance}
below, which will be the centerpiece of the KAM theorem of
Section~\ref{sec:KAMthm}.

\begin{defn} \label{def:inv}
Let $f_\param \colon {\mathcal M}\rightarrow{\mathcal M}$ be a family
of conformally symplectic maps.  A KAM attractor with frequency
$\omega$ is an invariant torus which is described by an embedding $K
\colon {\mathbb T}\rightarrow{\mathcal M}$ and a drift parameter
$\param$, which satisfy the following invariance equation for
$\theta\in{\mathbb T}$:
\begin{equation} \label{invariance}
f_\param \circ K(\theta) = K(\theta+\omega)\ .
\end{equation}
\end{defn}

We remark that solving equation \eqref{invariance} will require to
determine both $K$ and $\param$.

Denoting by $T_\omega$ the shift by $\omega$ such that for a function
$K$, we have $(K\circ T_\omega)(\theta)=K(\theta+\omega)$, then the
invariance equation \equ{invariance} can be written as
$$
f_\ecc\circ K=K\circ T_\omega\ .
$$

\subsection{The KAM theorem} \label{sec:KAMthm}

The KAM statement provided in \cite{CCL20} applies to two-dimensional
maps and, although it has been applied to the dissipative standard
map, the formulation of the KAM theorem was given for a general
system. Therefore, we can apply the main theorem stated in
\cite{CCL20} to the Poincar\'e map of the spin-orbit problem
\equ{eq.spin-diss}.

The KAM theorem in \cite{CCL20} gives explicit conditions that ensure
that, given an approximate solution, there is a true solution.  This
requires the computation of several constants that we list in
Appendix~\ref{app:constants} to make the paper self contained.  If the
map was given by an explict formula (as it was the case in
\cite{CCL20}) some of the constants can be obtained using calculus. In
our case, since the map is obtained integrating an ODE, we obtain the
estimates integrating the equation in a complex domain.

Having fixed a Diophantine frequency $\omega$ and after computing the
value of the conformal factor $\lambda$, we look for an embedding $K$
and a drift parameter $\param$ which satisfy the invariance equation
\equ{invariance}. The solution can be obtained under a non-degeneracy
condition (see \ref{main.H3} in Theorem~\ref{main}).

In the spin-orbit problem, the description of the computation of the
solution is given in Section~\ref{sec:initial}, while the verification
of the KAM conditions is provided in Section~\ref{sec:estimates}.

\vskip.1in

Let us assume that we start with an approximate solution
$(K_0,\param_0)$ which satisfies the invariance equation
\eqref{invariance} up to an error term $E_0$, that is,
\begin{equation}
\label{approx}
 E_0(\theta) = f_{\param_0} \circ K_0(\theta) - K_0(\theta+\omega)\ .
\end{equation}
Before stating the main theorem, we need to introduce the following
auxiliary quantities:
\begin{equation}
\label{def}
\begin{split}
N_0(\theta) &\coloneq (DK_0(\theta)^\top DK_0(\theta))^{-1}\ , \\
M_0(\theta) &\coloneq [ DK_0(\theta)\ |\  J^{-1}\circ K_0(\theta)\ DK_0(\theta) N_0(\theta)]\ , \\
S_0(\theta) & \coloneq ((DK_0 N_0)\circ T_\omega)^\top(\theta)
Df_{\param_0} \circ K_0(\theta) J^{-1}\circ
K_0(\theta)DK_0(\theta)N_0(\theta)\ ,
\end{split}
\end{equation}
where the superscript $\top$ denotes transposition and the matrix $J$
is the matrix representation of the symplectic form,
\[
\Omega_z(u,v) = \langle u, J(z) v \rangle,
\]
with $z\in \mathcal M$. For the applications we have in mind, $J$ is
constant and it is defined as
\begin{equation}
 \label{matrixJ}
 J =
 \begin{pmatrix}
  0 & 1 \\ -1 & 0
 \end{pmatrix}\ .
\end{equation}

Theorem~\ref{main} is a constructive version of Theorem 20 in
\cite{CallejaCL11} and it applies to mapping systems, like the
Poincar\'e map $P _\ecc$ defined in \eqref{eq.Pe} associated to
\eqref{eq.spin-diss}. In this case the conformal factor $\lambda$ only
depends on the dissipation $\dis$ and the eccentricity $\ecc$, and the
map $P _\ecc$ depends on the three parameter $\dis$, $\eps$, and
$\ecc$.

\begin{thm}\label{main}
Let $\Lambda$ be an open subset of $\mathbb{R}$ and for all $\param
\in \Lambda$, let $f_\param \colon {\mathcal M}\rightarrow {\mathcal
  M}$ be a conformally symplectic map defined on the manifold
${\mathcal M} = B\times {\mathbb T}$; here $B\subset {\mathbb R}$
denotes an open and simply connected domain with smooth boundary.
Assume that $f_\param$ is analytic on an open connected domain
$\mathcal{C}\subset{\mathbb C}\times{\mathbb C}/{\mathbb Z}$.  Assume
the following hypotheses.
\begin{enumerate}
\renewcommand*{\theenumi}{\bfseries{H\arabic{enumi}}}
\renewcommand*{\labelenumi}{\theenumi.}
 \item \label{main.H1} The frequency $\omega$ is Diophantine as in
   \equ{DC}, namely $\omega\in {\mathcal D}(\nu,\tau)$.
 \item \label{main.H2} The approximate solution $(K_0,\param_0)$,
   $K_0\in\mathcal{A}_{\rho_0}$ for some $\rho_0>0$ and
   $\param_0\in\Lambda$, satisfies \eqref{invariance} up to an error
   function $E_0=E_0(\theta)$ as in \eqref{approx}. We denote by
   $\varepsilon_0$ the size of the error function, that is,
   \[
    \varepsilon_0 \coloneq \|E_0\|_{\rho_0}\ .
   \]
 \item \label{main.H3} Assume that the following non--degeneracy
   condition is fulfilled:
\[
\det
\begin{pmatrix}
 {\overline S}_0 & {\overline {S_0(B_{b0})^0}}+\overline{\widetilde
    A_0^{(1)}} \\ \lambda-1 & \overline{\widetilde A_0^{(2)}}
\end{pmatrix} \ne 0\ ,
\]
 where $S_0$ is defined in \eqref{def}, $\widetilde A_0^{(1)}$,
 $\widetilde A_0^{(2)}$ are the first and second elements of
\[\widetilde A_0 = M_0^{-1}\circ T_{\omega} D_\param f_{\param_0} \circ
K_0\ ,\] $(B_{b0})^0$ is the solution (with zero average in the
$\lambda = 1$ case) of the equation \[\lambda
(B_{b0})^0-(B_{b0})^0\circ T_\omega=-(\widetilde A_0^{(2)})^0\ ,\] and
$(\widetilde A_0^{(2)})^0$ is the zero average part of $\widetilde
A_0^{(2)}$.

Then, let ${\mathcal T}_0$ be the \emph{twist constant} defined as
\[
{\mathcal T}_0  \coloneq \left \|
\begin{pmatrix}
  {\overline S}_0 & {\overline {S_0(B_{b0})^0}}+\overline{\widetilde
    A_0^{(1)}} \\ \lambda-1 & \overline{\widetilde A_0^{(2)}}
\end{pmatrix}^{-1} \right \|\ .
\]
 \item \label{main.H4} Assume that for some $\zeta>0$ we have
\[ {\rm dist}( \param_0, \partial \Lambda) \ge \zeta\ ,\qquad {\rm
  dist}(K_0({\mathbb T}_{\rho_0}), \partial\mathcal{C}) \ge \zeta\ .
\]
 \item \label{main.H5} Let $\delta_0$ be such that
   $0<\delta_0<\rho_0$.  Introduce the quantity $\kappa_\param
   \coloneq 4 C_{\sigma 0}$ with $C_{\sigma 0}$ constant (see
   Appendix~\ref{app:constants}). Define the quantities
\begin{equation}
 \label{eq.Qs}
 \begin{split}
Q_z&\coloneq \sup_{z\in\mathcal{C}}|Df_{\param_0}(z)|\ , \\
Q_{\param}&\coloneq \sup_{z\in\mathcal{C},\param\in\Lambda,|\param-\param_0|<2\kappa_\param\varepsilon_0}|D_\param f_{\param}(z)|\ , \\
Q_{zz}&\coloneq \sup_{z\in\mathcal{C}}|D^2 f_{\param_0}(z)|\ , \\
Q_{\param z}&\coloneq \sup_{z\in\mathcal{C}}|DD_\param f_{\param_0}(z)|\ , \\
Q_{zzz}&\coloneq \sup_{z\in\mathcal{C}}|D^3 f_{\param_0}(z)|\ ,  \\
Q_{\param zz}&\coloneq \sup_{z\in\mathcal{C},\param\in\Lambda,|\param-\param_0|<2\kappa_\param\varepsilon_0}|D^2 D_\param f_{\param}(z)|\ ,  \\
Q_{z\param }&\coloneq \sup_{z\in\mathcal{C},\param\in\Lambda,|\param-\param_0|<2\kappa_\param\varepsilon_0}|D_\param Df_{\param}(z)|\ , \\
Q_{\param\param}&\coloneq \sup_{z\in\mathcal{C},\param\in\Lambda,|\param-\param_0|<2\kappa_\param\varepsilon_0}|D^2_\param f_{\param}(z)|\ , \\
Q_{zz\param}&\coloneq \sup_{z\in\mathcal{C},\param\in\Lambda,|\param-\param_0|<2\kappa_\param\varepsilon_0}|D_\param D^2 f_{\param}(z)|\ , \\
Q_{\param\param z}&\coloneq \sup_{z\in\mathcal{C},\param\in\Lambda,|\param-\param_0|<2\kappa_\param\varepsilon_0}|DD_\param^2 f_{\param}(z)|\ , \\
Q_{\param\param\param}&\coloneq \sup_{z\in\mathcal{C},\param\in\Lambda,|\param-\param_0|<2\kappa_\param\varepsilon_0}|D_\param^3 f_{\param}(z)|\ ,\\
Q_{E 0}&\coloneq
{1\over2}\max\Big\{\|D^2E_0\|_{{\rho_0-\delta_0}},\|DD_\param
E_0\|_{{\rho_0-\delta_0}}, \|D^2_\param
E_0\|_{{\rho_0-\delta_0}}\Big\}\ .
 \end{split}
\end{equation}

Assume that $\varepsilon_0$ is such that the following smallness
conditions are satisfied for real constants $C_{\eta 0}$,
$C_{{\mathcal E} 0}$, $C_{d0}$, $C_{\sigma 0}$, $C_\sigma$,
$C_{W0}$, $C_W$ and $C_{\mathcal{R}}$ (see
Appendix~\ref{app:constants}):
\begin{align}
 C_{\eta 0}\,\nu^{-1}\delta_0^{-\tau}\varepsilon_0 &< \zeta\ , \label{C1} \\
2^{3\tau+4}\,C_{{\mathcal E} 0}\ \nu^{-2}\
\delta_0^{-2\tau} \varepsilon_0 &\leq 1\ , \label{C2} \\
4C_{d0} \nu^{-1}\delta_0^{-\tau} \varepsilon_0 &< \zeta\ ,  \label{C3} \\
4C_{\sigma 0}\varepsilon_0 &< \zeta\ ,   \label{C4} \\
\|N_0\|_{\rho_0}\ (2\|DK_0\|_{\rho_0}+D_K)\ D_K &< 1 \ , \label{condbT} \\
4 Q_{z\param 0} C_{\sigma 0} \varepsilon_0 &< Q_z\ ,  \label{Cnew1} \\
4 Q_{\param\param} C_{\sigma 0} \varepsilon_0 &< Q_{\param }\ ,   \label{Cnew2} \\
C_\sigma\ D_K &\leq C_{\sigma 0}\ , \label{C8} \\
D_K(C_{W0}+\|M_0\|_{\rho_0}C_W+C_W D_K) &\leq C_{d0}\ , \label{C9} \\
D_K\ \Big(C_W\ \nu \delta_0^{-1+\tau}+C_{\mathcal{R}}\Big) &\leq
C_{{\mathcal E} 0}\ ,  \label{C10}
\end{align}
\noindent
where $D_K$ is given by
\begin{equation} \label{DK1} D_K\coloneq 4C_{d0}\
\nu^{-1}\delta_0^{-\tau-1}\ \varepsilon_0\ .
\end{equation}
\end{enumerate}
Then, there exists an exact solution $(K_*,\param_*)$ of
\equ{invariance} satisfying
\[
f_{\param_*} \circ K_* - K_* \circ T_\omega = 0\ .
\]
The following inequalities show that the quantities $(K_*, \param_*)$
are close to $(K_0, \param_0)$:
\begin{equation}
 \label{kemu}
 \begin{split}
  \| K_* - K_0 \|_{\rho_0-\delta_0} &\leq 4 C_{d0} \nu^{-1}
  \delta_0^{-\tau} \|E_0\|_{\rho_0}\ , \\ | \param_* - \param_0| &\leq
  4 C_{\sigma 0} \|E_0\|_{\rho_0}\ ,
 \end{split}
\end{equation}
where $C _{d0}$ and $C _{\sigma 0}$ are given explicitly in
Appendix~\ref{app:constants}.
\end{thm}

For simplicity of exposition, we report the explicit expressions of
the constants entering Theorem~\ref{main} in
Appendix~\ref{app:constants}. They are obtained making a constructive
version of the KAM proof given in \cite{CallejaCL11}.  We refer to
\cite{CCL20} for the proof of Theorem~\ref{main}.


\subsection{A sketch of the proof of Theorem~\ref{main}}
\label{sec:sketch}

We present a sketch of the proof of Theorem~\ref{main} that we split
into five main steps, all of them giving explicit estimates of the
quantities involved. Although we do not enter into the details of the
proof, which is quite long and technical (see \cite{CCL20}), we
provide an overview of the proof which motivates the assumptions
\ref{main.H1}-\ref{main.H5} as well as the smallness conditions
\equ{C1}-\equ{C10}.

\subsubsection{Step 1: the approximate solution.}

We denote by $(K,e)$ an embedding function and a drift term satisfying
approximately the invariance equation with an error term $E$:
\begin{equation} \label{invE}
f_e\circ K(\theta)-K(\theta+\omega)=E(\theta)\ .
\end{equation}
All one-dimensional tori are Lagrangian invariant tori, namely they
satisfy $K^*\Omega=0$, which in coordinates is given by
$$
DK^T(\theta)\ J\circ K(\theta)\ DK(\theta)=0\ .
$$
This expression implies that the tangent space can be decomposed as
the sum of the range of $DK(\theta)$ and the range of $V(\theta)$,
where $V$ is given by
$$
V(\theta)=J^{-1}\circ K(\theta)\ DK(\theta) N(\theta)
$$
with $N(\theta)=(DK(\theta)^\top DK(\theta))^{-1}$.

Next, we define the quantity $M$ as a juxtaposition of $DK$ and $V$,
i.e.,
\begin{equation} \label{Mdef}
M(\theta)=[DK(\theta)\ |\  V(\theta)]\ .
\end{equation}
Then, it can be shown that, up to a remainder $R$, the action of the
derivative of the map over $M$ is just a shift of $M$ multiplied by a
matrix. Precisely, one can prove that (\cite{CCL20}):
\begin{equation} \label{R}
Df_e \circ K(\theta)\ M(\theta) = M(\theta + \omega)
\left(%
\begin{array}{cc}
  {\rm Id} & S(\theta)\\
  0 & \lambda {\rm Id} \\
\end{array}%
\right) +R(\theta)\ .
\end{equation}
This result will be used in Step 2 to reduce \equ{invE} to a constant
coefficient equation, that will be solved under assumptions
\ref{main.H1} and \ref{main.H3}.

\subsubsection{Step 2: a new approximation.}

Starting from the initial approximation $(K,e)$, we introduce a new
approximation $(K',e')$ defined adding to $(K,e)$ some corrections
$(W,\sigma)$ as $K'=K+M W$ and $e'=e+\sigma$.  We denote by $E'$ the
error function associated to $(K',e')$, satisfying the equation:
\begin{equation} \label{APPR-INVp}
f_{e'}\circ K'(\theta)-K'(\theta+\omega)=E'(\theta)\ .
\end{equation}
Next, we proceed to expand \equ{APPR-INVp} in Taylor series, which gives:
\beqano
&&f_e\circ K(\theta)+Df_e\circ K(\theta)\ M(\theta)W(\theta)+D_{e}  f_{e}
\circ K(\theta){\sigma}\nonumber\\
&&\qquad -K(\theta+\omega)-
M(\theta+\omega)\ W(\theta+\omega)+h.o.t.=E'(\theta)\ .
\eeqano
Using \equ{invE}, we can guarantee that $E'$ is quadratically
smaller provided that the following relation is satisfied:
\begin{equation} \label{quad}
Df_e \circ K(\theta)\ M(\theta)W(\theta) - M(\theta+\omega)\ W(\theta+\omega) +
D_e  f_e \circ K(\theta)\sigma = - E(\theta)\ .
\end{equation}
We remark that condition \equ{C1} provides an estimate of the error
$E'$ associated to $(K',e')$.

Using \equ{quad} and \equ{R}, we obtain that
$$
Df_{e} \circ K(\theta)\ M(\theta) = M(\theta +\omega)
\left(%
\begin{array}{cc}
  {\rm Id} & S(\theta) \\
  0 & \lambda {\rm Id} \\
\end{array}%
\right)+R(\theta)\ ,
$$
which provides the following equations for $W$ and $e$:
\begin{equation} \label{We}
M(\theta +\omega)
\left(%
\begin{array}{cc}
  {\rm Id} & S(\theta)\\
  0 & \lambda {\rm Id} \\
\end{array}%
\right) \ W(\theta) - M(\theta+\omega)\ W(\theta+\omega) = - E(\theta)- D_{e}
f_{e} \circ K(\theta){\sigma}\ .
\end{equation}
Next, we multiply by $M(\theta+\omega)^{-1}$ and write \equ{We}
for the components $W_1$, $W_2$,
$\tilde E_1$, $\tilde E_2$,
$\tilde A_1$, $\tilde A_2$,
 of $W$, $\tilde E$, and $\tilde A$ as
\begin{equation} \label{W12}
\left(%
\begin{array}{cc}
  {\rm Id} & S(\theta) \\
  0 & \lambda {\rm Id} \\
\end{array}%
\right) \left(\begin{array}{c}
  W_1(\theta) \\
  W_2(\theta) \\
\end{array}\right)-
\left(\begin{array}{c}
  W_1(\theta+\omega) \\
  W_2(\theta+\omega) \\
\end{array}\right)=
\left(\begin{array}{c}
  -\tilde E_1(\theta)-\tilde A_1(\theta)\sigma \\
  -\tilde E_2(\theta)-\tilde A_2(\theta)\sigma \\
\end{array}\right)\ ,
\end{equation}
where we define
\begin{align*}
 \tilde E_j(\theta) &= -(M(\theta+\omega)^{-1}E)_j\nonumber \\
\tilde A_j(\theta) &= (M(\theta +\omega)^{-1}D_e f_e\circ K)_j\ ,
\end{align*}
for $j =1,2$.
We now make explicit \equ{W12} for the components $W_1$, $W_2$ and
$\sigma$, so to obtain the following cohomological equations:
\begin{equation}
\label{B}
 \begin{split}
  W_1(\theta)-W_1(\theta+\omega) &= -\widetilde E_1(\theta)-S(\theta) W_2(\theta)-\widetilde A_1(\theta)\,
\sigma \\
\lambda W_2(\theta)-W_2(\theta+\omega) &= -\widetilde E_2(\theta)-\widetilde A_2(\theta)\, \sigma\ .
 \end{split}
\end{equation}

\subsubsection{Step 3: Determination of the new approximate solution.}

The solution of equations \equ{B} allow us to determine the unknowns
$W_1$, $W_2$ and $\sigma$ that give the corrections to determine the
new approximate solution.

To solve the first equation of \equ{B}, we use assumption
\ref{main.H1} on the Diophantine property of the frequency and
assumption \ref{main.H3}, expressing the non-degeneracy that allows us
to solve the linear system \equ{system} below. The second equation of
\equ{B} can instead be solved by an elementary contraction mapping
argument for any $|\lambda|\ne 1$ and for all real frequencies.

Let us write $W_2$ as $W_2 = \langle W_2\rangle + B^0 + \tilde B^0
\sigma$.  Taking the average of both equations \equ{B}, we obtain the
equations
\begin{equation} \label{system}
\begin{pmatrix}
 \langle S\rangle & {\langle {SB^0}}\rangle +\langle {\widetilde A_1}\rangle  \\
  (\lambda-1){\rm Id} & \langle {\widetilde A_2}\rangle
\end{pmatrix}
\begin{pmatrix}
 \langle W_2\rangle  \\
\sigma 
\end{pmatrix} = 
\begin{pmatrix}
 -\langle S \tilde B^0 \rangle- \langle \widetilde E_1\rangle \\
-\langle \widetilde {E_2}\rangle
\end{pmatrix} \ ,
\end{equation}
which can be solved to give $\langle W_2\rangle$ and $\sigma$ under
the non-degeneracy condition \ref{main.H3}.

Once the solution of \equ{system} is obtained, we proceed to solve the
second of \equ{B} to determine $W_2$; such equation can be solved for
any $|\lambda|\ne 1$ by a contraction mapping argument.

Then, we proceed to solve the first equation of \equ{B} for $W_1$:
since it involves small divisors, we can solve the equation under the
Diophantine assumption \ref{main.H1}. The quantities
$\|W_1\|_{\rho-\delta}$ and $\|W_2\|_{\rho-\delta}$ can be bounded by
$\|E\|_\rho$ by using Cauchy estimates for the cohomological equations
\equ{B}.

The error $E'$ associated to the new solution can be bounded on a
domain of size $\rho-\delta$ by the square of the error $E$ on the
domain of size $\rho$ as
\[
\|E'\|_{\rho-\delta}\leq C_E\delta^{-2\tau}\|E\|_\rho^2\ ,\qquad C_E>0\ ,
\]
showing that the new error of the procedure is quadratic in the
original error. Assumption \ref{main.H4} is needed to obtain such a
bound.

\subsubsection{Step 4: iteration and convergence.}

We proceed to iterate the procedure presented in Step 3 to obtain a
sequence of new solutions, say $\{K_j,e_j\}$, and their associated
invarance equation error, say $E _j$. We prove that the errors tends
to zero (in suitable norms) as $j \to \infty$ and thus the solution
sequence converges to the true solution. The proof consists in
implementing an abstract implicit function theorem, alternating the
iteration with carefully chosen smoothing operators for analytic
functions. The smoothing is obtained by rescaling domains where the
functions are defined at each step.  In particular, we can define as
$\rho_j$ the size of the analyticity domain associated to the solution
$\{K_j,e_j\}$ by introducing a shrinking parameter $\delta_j$ and
setting
\[
\rho_0=\rho\ ,\qquad \delta_j={\rho_0\over {2^{j+2}}}\ ,\qquad
\rho_{j+1}=\rho_j-\delta_j\ ,\qquad j\geq 0\ .
\]
Then, we can show that for $a,b>0$ and $C_E'>0$, we have
\[
\|E _{j+1}\|_{\rho_{j+1}}\leq C_E'\ \nu^a\delta_j^b\
\|E _j\|_{\rho_j}^2\ .
\]
If the quantity $\varepsilon _0\equiv \|E _0\|_{\rho_0}$ is
sufficiently small, then we conclude that
\begin{equation} \label{Ke}
\|K_j-K_0\|_{\rho_j}\leq C_K\varepsilon_0\ , \qquad |e_j-e_0|\leq
C_\mu\varepsilon_0
\end{equation}
for some constants $C_K,C_\mu>0$. The inequalities \equ{C1}-\equ{C10}
of Theorem~\ref{main} allow to obtain \equ{Ke} as well as to ensure
that the procedure can be iterated and that it converges to the true
solution.

\subsubsection{Step 5: local uniqueness.}

Under smallness conditions, one can prove that, if there exist two
solutions $(K_a,e_a)$, $(K_b,e_b)$, then there exists $\psi\in\R$ such
that
\[
K_b(\theta)=K_a(\theta+\psi)\qquad {\rm and}\qquad e_a=e_b\ .
\]

\subsection{The algorithm and the initial invariant curve}
\label{sec:initial}

Theorem~\ref{main} provides an explicit algorithm working as follows:
for a fixed frequency $\omg$ and from an approximate solution
$(K_0,e_0)$ satisfying the invariance equation with error term $E_0$,
one can construct a new approximation $(K_1,e_1)$ satisfying the
invariance equation with a new error term $E_1$ which is quadratically
smaller than $E_0$, just taking derivatives and performing algebraic
operations.  The new approximation is obtained by solving suitable
cohomological equations, under the non-degeneracy condition
\ref{main.H3}. The algorithm is presented in detail in \cite{CCGL20a}
for the spin-orbit problem and it is recalled in
Appendix~\ref{app:algorithm}.

\vskip.1in

In the following, we will consider two frequencies defined as
\begin{align}
 \label{eq.omg1}
  \omega_1 &\coloneq \gamma _g^+
  \intertext{and}
 \label{eq.omg2}
  \omega_2 &\coloneq 1 + \frac{1}{2 + \gamma _g^-} \ ,
\end{align}
where $\gamma _g^\pm \coloneq \frac{\sqrt{5}\pm 1}{2}$.  Both
frequencies are Diophantine, in the sense of
Definition~\ref{defn.Diophantine}, with constant $\nu =
(\frac{3-\sqrt{5}}{2})^{-1}$ and exponent $\tau = 1$.

The application of Theorem~\ref{main} consists in the steps given
below.

\vskip.1in

\begin{enumerate}
\renewcommand{\theenumi}{(\roman{enumi})}
\renewcommand{\labelenumi}{\theenumi}
\renewcommand{\itemsep}{.1in}
 \item We fix the Diophantine frequency as one of the choices in
   \equ{eq.omg1} or \equ{eq.omg2}.
 \item \label{stepii} We provide the initial values $K _0$ and $\ecc
   _0$, selecting the eccentricity and the initial condition as
   follows.  First, we select the eccentricity by choosing the value
   that corresponds to the fixed frequency. This is achieved by
   integrating equation \eqref{eq.spin-diss} with an initial guess of
   $e_0$ and initial conditions $x(0)=0$ and we fix $y(0)={\bar
     N}(\ecc_0)/{\bar L}(\ecc_0)$, which is the value that we obtain
   when the dissipation disappears in the averaged model
   \equ{eq.avg-spinxy}. After a transient time $t$ (so that the system
   evolves on the attractor), we compute the frequency over $N_{it}$
   additional iterations as $\omega={1\over N_{it}}
   \sum_{j=1}^{N_{it}} y(t+2\pi j)$. Once the approximated initial
   eccentricity for the desired frequency has been obtained, we
   iterate the Poincar\'e map (after another suitable transient) and
   we obtain the initial approximation of the invariant curve by
   fitting the discrete points.
 \item \label{stepiii} We iterate Algorithm~\ref{alg.newton} to obtain
   a more accurate approximation $(K_a,e_a)$ satisfying the invariance
   equation with an error whose norm is sufficiently small.
 \item \label{stepiv} We compute the norms of the quantities appearing
   in Theorem~\ref{main} and detailed in
   Appendix~\ref{app:quantities1} for $\omega_1$ and
   Appendix~\ref{app:quantities2} for $\omega_2$.
 \item We check the conditions \equ{C1}-\equ{C10} in
   Theorem~\ref{main}.  If they are satisfied, we conclude the
   procedure, otherwise we change some of the parameters (e.g., $\rho$
   and $\delta$) and we try to optimize the final result.
\end{enumerate}

Further details of the steps \ref{stepii} and \ref{stepiii} can be
found in \cite{CCGL20a}, which contains also the computation of the
variational equations with respect to the initial conditions and the
parameter $\ecc$ in \eqref{eq.spin-diss}. The rotation number in
\ref{stepii} can be computed more efficiently (with smaller $N _{it}$)
by \cite{Sander2017}. The variational equations are needed in step
\ref{alg.newton-PSA} of the Algorithm~\ref{alg.newton} in
Appendix~\ref{app:algorithm} as well as for some of the quantities in
\ref{stepiv}.


\subsection{Continuation method}
\label{sec.conti} 

Algorithm~\ref{alg.newton} can be used as a corrector for a
continuation method of the invariant torus and its drift. In the
spin-orbit problem, we use the eccentricity $\ecc$ as the adjustable
parameter required by the quasi-Newton method and the perturbative
parameter $\eps$ in \eqref{eq.spin-diss} as the continuation
parameter.
The continuation consists in increasing $\eps$ by a stepsize, say
$\eps _h$, and run the Algorithm~\ref{alg.newton} again with a given
Newton's tolerance $\tilde\epsilon$. Thus at each continuation step,
it succeed, we obtain a new embedding of the torus and a new corrected
eccentricity.

If $\eps + \eps _h$ converges, we increase $\eps _h$ for the next
continuation step. Otherwise, we do not accept $\eps + \eps _h$ as a
solution, we decrease $\eps _h$, and we use Algorithm~\ref{alg.newton}
with the new value of $\eps + \eps _h$.  In both cases we perform a
Lagrange interpolation of the previous two or three steps in order to
provide a better initial guess of $K$ and $\ecc$ for the next
iteration.

In all the above process, a refinement of the grid in the coordinate
$\theta$ may be required. In our implementation, we consider necessary
to increase the number of Fourier coefficients when some of the
following two cases arise.

The first one is when the accuracy tests, detailed in Section 5.3 of
\cite{CCGL20a}, fail.  In short, the accuracy tests are aimed to
control different sources of error, precisely:

\vskip.1in

\begin{enumerate}
\renewcommand{\itemsep}{.1in}
 \item the error of the invariance equation on a table of values;
 \item the error in the numerical integration, for which we introduce
   absolute and relative tolerances;
 \item the error in the grid over the coordinate $\theta$, which is
   controlled by checking the last coefficients of the truncated
   Fourier series as well as the Sobolev norm of the tail;
 \item the interpolation error, which is controlled by providing an
   estimate of it and by changing the size of the grid, when the error
   becomes too large.
\end{enumerate}

\vskip.1in

The second situation is when the continuation step fails consecutively
two times which may require to decrease the stepsize $\eps_h$,
especially when we are getting close to the breakdown. Thus, the
continuation procedure will stop when the maximum number of remeshing
is reached, in our case $2^{14}$ Fourier modes.

\begin{figure}[ht]
 \resizebox{\columnwidth}{!}{\input{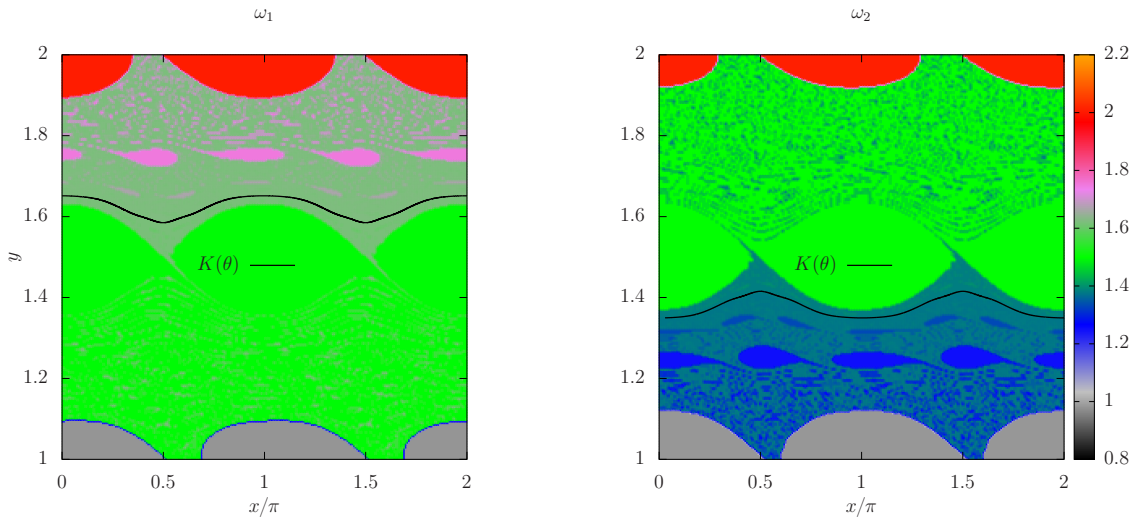}}
 \caption{Basins of rotation number given by color-scale for the
   parameters in \eqref{eq.conti-example-omg2} (left) and
   \eqref{eq.conti-example-omg1} (right) joined with the of
   Algorithm~\ref{alg.newton} that show the invariant attractor (in
   black) after a continuation of $\eps$ starting with $\eps=10^{-4}$
   and $\ecc=0.3150628$ (left) and $\ecc=0.2502068$ (right).}
 \label{fig.conti-spin}
\end{figure}

Figure~\ref{fig.conti-spin} displays the results of the KAM torus
(black curve) after a continuation starting at $\eps = 10^{-4}$ and a
fixed dissipation $\dis=10^{-3}$. The computation has been done with a
multi-precision arithmetic with $170$ bits, i.e. around $50$ digits of
accuracy, a Newton's tolerance of $\tilde\epsilon=10^{-35}$, and a
parallelization of the integration of the Poincar\'e map as detailed
in Section 5.5 of \cite{CCGL20a}.

We emphasize that we checked the final result by changing the number
of digits of accuracy; in other words, keeping the same Newton's
tolerance $\tilde\epsilon$ and the same integration's tolerance, we
have performed the last continuation step, checking that it is
satisfied with $50$, $55$, and $60$ digits of accuracy.

The values of the Fourier modes $n _\theta$, the dissipation $\dis$,
the eccentricity $\ecc$, and the perturbing parameter $\eps$ are
reported below for the Diophantine frequencies $\omg _1$ in
\eqref{eq.omg1} and $\omg _2$ in \eqref{eq.omg2}. \newline For $\omg
_1$, the last successful Newton continuation step was reached in less
than 3 min using 32 CPUs with final values:
\begin{equation}
\label{eq.conti-example-omg1}
 \begin{split}
  n _\theta &= 16384\ , \\
  \dis &= 10^{-3}\ , \\
  \ecc &= \mathtt{0.31675286891174832107186084513865661761571784973618}\ , \\
  \eps &= \mathtt{0.011632963641877116367716112642948530559675531382297}\ .
 \end{split}
\end{equation}
For $\omg _2$ we got the last successful torus in less than 5 min
using 35 CPUs and with values:
\begin{equation}
\label{eq.conti-example-omg2}
 \begin{split}
  n _\theta &= 4096\ , \\
  \dis &= 10^{-3}\ , \\
  \ecc &= \mathtt{0.24824740823563165902227100091869770425731996450084}\ , \\
  \eps &= \mathtt{0.012697630024415883032123830013667613509009950826168}\ .
 \end{split}
\end{equation}
Figure~\ref{fig.conti-spin} provides also the basins of the rotation
numbers, namely the frequency given through a color scale for
different initial conditions $(x_0,y_0)$. In particular, we take a
grid of $500\times 500$ initial conditions within the window
$[0,2\pi)\times [1,2]$ and we compute the frequency as described in
  step \ref{stepii} of Section~\ref{sec:initial}. We remark that the
  computation of the frequency has been optimized using the method
  described in \cite{Sander2017} which is implemented and detailed for
  the spin-orbit case in the companion paper \cite{CCGL20c}.

\section{Estimates on the $Q$ quantities of the KAM theorem~\ref{main}}
\label{sec:Qestimates}


The main difference in the explicit derivation of the KAM estimates
presented in \cite{CCL20} between the standard map and the spin-orbit
problem is the computation of the $Q$ constants defined in \equ{eq.Qs}
of Theorem~\ref{main}.  Almost all of them are zero for the standard
map, while for the spin-orbit problem we need to compute them as
detailed in Sections~\ref{sec:QE0} and \ref{sec:Q} below.

It is also important to describe carefully the boundary of the domain
$\mathcal{C}$ and, in particular, the value $\zeta$ in \ref{main.H4}
which is needed for the inequalities \eqref{C1}--\eqref{C10}.

\subsection{The computation of $Q _{E0}$}
\label{sec:QE0}
We need to give a bound of the quantity
\begin{equation}
\label{eq.QE0}
 Q _{E _0} \coloneq \frac{1}{2} \max \bigl\{ \| D^2 E _0 \| _{\rho _0 -
   \delta _0}, \| D _ \ecc D E _0 \| _{\rho _0 - \delta _0}, \| D^2
 _\ecc E _0 \| _{\rho _0 - \delta _0} \bigr\}\ ,
\end{equation}
where $E _0$ is defined in terms of the numerical approximated
solution $(K _0, \ecc _0)$ of Theorem~\ref{main}. In the case of the
spin-orbit problem, $E _0$ is given by
\begin{align}
 \mathcal E(\theta) &\coloneq (\Psi^{-1} _{\ecc _0}\circ G _{\ecc _0})^{1} (
 \Psi _{\ecc _0}\circ K _0(\theta)) - K _0^1 (\theta +
 \omega) \ , \nonumber \\
 \label{eq.E01}
 E _0^1(\theta) &\coloneq \mathcal{E}(\theta) - \lfloor \mathcal{E}(\theta) +
 0.5 \rfloor \ , \\
 \label{eq.E02}
 E _0^2(\theta) &\coloneq (\Psi^{-1} _{\ecc _0}\circ G _{\ecc _0})^{2} ( \Psi
 _{\ecc _0}\circ K _0(\theta)) - K _0^2 (\theta + \omega)\ ,
\end{align}
where $\lfloor \, \cdot \, \rfloor$ denotes the floor function, $\ecc
_0$ is the eccentricity value, $G _{\ecc _0} = \Psi _{\ecc _0} \circ P
_{\ecc _0} \circ \Psi _{\ecc _0}^{-1}$ with $P _{\ecc _0}$ being the
$2\pi$-time flow of \eqref{eq.spin-diss} and $\Psi _{\ecc _0}$ given
in \eqref{eq.xy2bg}. The superscripts ${}^{1}$ and ${}^{2}$ mean the
components of the vectors in $\mathbb{R}^2$. Note that the floor
function in \eqref{eq.E01} is needed since $x(t)$ in
\eqref{eq.spin-diss} is given modulus $2\pi$, that due to $\Psi _{\ecc
  _0}$, in fact, it is modulus $1$. Therefore $E _0^1$ gives values in
$[-1/2,1/2]$.

To compute $D ^2 E _0$ we can either differentiate the Fourier series
with respect to $\theta$ or to use jet transport, which, roughly
speaking, means to overload the numerical integrator with a
multivariate polynomial manipulator.  We are going to use the jet
transport because we also need to get the variation with respect to
the eccentricity, i.e., $D _\ecc D E _0$ and $D _\ecc E _0$. In order
to get the quantities automatically, we use jets\footnote{We follow
  the convention that a jet is encoded by the Taylor's coefficients at
  0.} of $2$ symbols, say $(s _1, s _2)$, and up to degree $2$, see
Appendix~\ref{sec.multipoly2}.  Indeed, for each $\theta$ in a mesh of
$\mathbb{T}$, we compute the flow given by
\begin{equation}
\label{eq.jet-spin}
 \Psi^{-1} _{\ecc _0 + s _2} \circ G _{\ecc _0 + s _2} \circ \Psi _{\ecc _0 + s
   _2}( K _0(\theta + s _1)) \ ,
\end{equation}
where
\begin{equation*}
 K _0(\theta + s _1) = K _0(\theta) + \partial _\theta K _0(\theta) s
 _1 + \tfrac{1}{2}\partial _\theta^2 K _0(\theta) s _1^2 \ .
\end{equation*}


\begin{rem}
 Jet transport will provide the normalized derivative of
 \eqref{eq.jet-spin}, so the $1/2$ in \eqref{eq.QE0} is automatically
 included in the coefficients of degree $2$ of
 \eqref{eq.jet-spin}. Notice that here, we can use the ad hoc
 polynomial manipulator described in Appendix~\ref{sec.multipoly2}.
\end{rem}

\begin{rem}
 The term $\lfloor\, \cdot + 0.5 \, \rfloor$ in \eqref{eq.E01} refers
 to the $\mathtt{round}$ function, namely the function that returns
 the nearest integer, but round halfway cases away from zero,
 regardless of the current rounding direction, and instead of the
 nearest integer in the $\mathtt{rint}$ function.

 Note that $\mathtt{round}$ has zero derivative except in
 $(\frac{1}{2}\mathbb{Z}) \setminus \{0\}$, where the derivative is
 not well-defined. However, we will consider (numerically) derivative
 zero also in these discontinuity points.
\end{rem}

\begin{rem}
 About the 2\textsuperscript{nd} derivatives for the term $K _0
 (\theta + \omega)$, the ones with respect to $\ecc$ are zero and the
 computation of $\partial _\theta ^2 K _0(\theta + \omega)$ is
 straightforward in the Fourier representation.
\end{rem}

\begin{rem}
 The computation of \eqref{eq.jet-spin} is fully parallelizable for
 each of the different values of $\theta$, which gives us a clear
 speed-up in the performance. Specially when the quantity is computed
 near to the breakdown parameter value that, generically, requires
 more Fourier modes.
\end{rem}

\subsection{The computation of the complex $Q$'s in Theorem~\ref{main}}
\label{sec:Q}

The quantities in the hypothesis \ref{main.H5} of Theorem \ref{main}
require to perform the integration of complex numbers, since the
initial conditions are in the complex domain $\mathcal{C}$, in fact,
in its boundary. The complexification of the spin-orbit model leads to
the \emph{complex spin-orbit} problem, see Section~\ref{sec:cspin},
which is given as a real 4-dimensional ODE system. This system
describes the evolution in time of the real and imaginary parts of
each of the variables in \eqref{eq.spin-diss}.

To address some of the freedoms in Theorem~\ref{main}, we devote our
attention in Section~\ref{sec:boundaryC} to provide a definition of a
possible domain $\mathcal{C}$ such that we can fulfill the hypothesis
\ref{main.H4}. The strategy will be to move this original freedom on
$\mathcal{C}$ to two new parameters, $\xi $ and $\alpha$, which are
going to be easier to handle.

Finally, we detail in Section~\ref{sec:stepsQ} the different steps to
approximate the $Q$ quantities of \ref{main.H5}.

\subsubsection{Complex spin-orbit problem} \label{sec:cspin}
The $Q$ quantities in \eqref{eq.Qs} are considered over the complex
domain $\mathcal{C}$ of Theorem~\ref{main}. This implies the need of
the complexification of the spin-orbit problem \eqref{eq.spin-diss},
which leads to a new system called the \emph{complex spin-orbit
  problem} given by the real ODE system
\begin{equation}
\label{eq.cspinxy}
 \begin{split}
  \frac{d}{dt} \re x(t) &= \re y(t) \ , \\
  \frac{d}{dt} \im x(t) &= \im y(t) \ , \\
  \frac{d}{dt} \re y(t) &= - \eps \Big(\frac{a}{r(t)}\Big)^3 \sin \bigl(2\re x(t) - 2f(t) \bigr) \cosh(2\im x(t)) - \dis
  \Big(\frac{a}{r(t)}\Big)^5 \bigl(\re y(t) - \frac{d}{dt}f(t) \bigr) \ ,
  \\
  \frac{d}{dt} \im y(t) &= - \eps \Big(\frac{a}{r(t)}\Big)^3 \cos \bigl(2\re x(t) - 2f(t) \bigr) \sinh(2\im x(t)) - \dis
  \Big(\frac{a}{r(t)}\Big)^5 \im y (t) \ ,
 \end{split}
\end{equation}
where \eqref{eq.cspinxy} has been deduced by taking the complex
numbers $x = \re x + \I \im y$ and $y = \re x + \I \im y$ in
\eqref{eq.spin-diss}. To obtain the above equations, we use the
relation
\[
 \sin (\re \alpha + \I \im \alpha) = \sin \re \alpha \cosh \im \alpha
 + \I \cos \re \alpha \sinh \im \alpha \ .
\]
Similarly to the real spin-orbit problem, see \cite{CCGL20a}, we
consider the temporal change of coordinates $t = u - \ecc \sin u$ to
make $u$ the independent variable, i.e.,
\begin{equation}
 \label{eq.cxy2bg}
 \begin{aligned}
  \re x (u - \ecc \sin u) &\eqcolon \re \beta(u)\ , & \re y (u - \ecc
  \sin u) &\eqcolon \re \gamma(u)/(1-\ecc \cos u)\ , \\ \im x (u -
  \ecc \sin u) &\eqcolon \im \beta(u)\ , & \im y (u - \ecc \sin u)
  &\eqcolon \im \gamma(u)/(1-\ecc \cos u)\ .
 \end{aligned}
\end{equation}
Thus, if $\widehat G _\ecc$ is the $2\pi$-time flow of the complex
spin-orbit problem with the coordinates $(\re \beta, \im \beta, \re
\gamma, \im \gamma)$, then we can recover the normalized $2\pi$-time
flow $\widehat P _\ecc$ of \eqref{eq.cspinxy} by the conjugacy given
by
\begin{equation*}
 \widehat \Psi _\ecc \coloneq
 2 \pi
 \begin{pmatrix}
  1 & 0 & 0 & 0 \\
  0 & 1 & 0 & 0 \\
  0 & 0 & 1 - \ecc & 0  \\
  0 & 0 & 0 & 1 - \ecc \\
 \end{pmatrix} \ .
\end{equation*}
Explicitly, we obtain:
\begin{equation}
 \label{eq.cPe}
 \widehat P _\ecc \coloneq \widehat \Psi ^{-1}_\ecc \circ \widehat G
 _\ecc \circ \widehat \Psi _\ecc \ .
\end{equation}
Therefore to get the different high variational flows involved in the
$Q$ quantities of \ref{main.H5}, we can use the jet transport
technique, see \cite{CCGL20a}, with jets of 5 symbols and up to order
3.

\subsubsection{Definition of the boundary of the complex domain $\mathcal{C}$} \label{sec:boundaryC}
The $Q$ quantities of the hypothesis \ref{main.H5} depend on the
boundary $\partial \mathcal{C}$, because the ODE \eqref{eq.spin-diss}
as well as \eqref{eq.cspinxy} are analytic. The only restriction on
this set $\partial \mathcal{C}$ is given in \ref{main.H4} which
relates the distance of the set
\begin{equation*}
 K _0(\mathbb{T} _{\rho _0}) \coloneq
 \biggl\{
 \begin{pmatrix}
  \theta + \I \sigma \\ 0
 \end{pmatrix} + \overline{K} _0(\theta + \I \sigma)
 \colon \theta \in \mathbb{T} \text{ and } |\sigma | \leq \rho _0 \biggr\}
\end{equation*}
with $\overline{K} _0$ denoting the periodic part of the mapping $K
_0$ which is continuously extented to the boundary of the set
$\mathbb{T} _{\rho _0}$ defined in \eqref{eq.Trho}.

Recall that the distance between sets is defined by
\begin{equation*}
 {\rm dist}(K _0(\mathbb{T} _{\rho _0}), \partial \mathcal{C})
 \coloneq \inf \{ d (x,y)\colon x \in K _0 (\mathbb{T}_{\rho _0})
 \text{ and } y \in \partial \mathcal{C}\} \ .
\end{equation*}

Hence, we consider $\mathcal{C}$ given in terms of a real region $\Xi$
in the plane and a real value $\alpha > 0$, as
\begin{equation*}
 \mathcal{C} \coloneq \{ (z _1, z _2) \in \mathbb{C}/\mathbb{Z} \times
 \mathbb{C}\colon \Re(z _1, z _2) \in \Xi, \, |\Im z _1 | \leq \alpha,
 \, |\Im z _2 | \leq \alpha\} \ .
\end{equation*}
The region $\Xi$ is bounded and we assume to be of the form
\begin{equation*}
 \Xi \coloneq \{ (\theta,\sigma ) \in \mathbb{T} \times \mathbb{R}
 \colon \psi _-(\theta) \leq \sigma \leq \psi _+(\theta) \}
\end{equation*}
for some real curves $\psi _-$ and $\psi _+$ such that
\begin{equation*}
 \psi _- \circ K_0^1(\theta) \leq K _0^2(\theta) \leq \psi _+ \circ
 K_0^1(\theta) \quad \text{for all } \theta \in \mathbb{T} \ .
\end{equation*}
For instance, fixed $\xi >0$, one can try to find $\psi _\pm$ solving
\begin{equation*}
 \psi _\pm \circ K_0^1(\theta) = K _0^2(\theta) \pm \xi \quad
 \text{for all } \theta \in \mathbb{T} \ .
\end{equation*}
Then, $K _0(\theta) \in \Xi$ for all $\theta$ in
$\mathbb{T}$. Heuristically, $K_0^1(\theta) = \theta +
\overline{K}_0^1(\theta) \approx \theta$, if $\overline{K}_0^1$ is
small and the composition by $K_0^1$ may be neglected. In fact, if we
allow constant values for $\psi _\pm$, we can just consider
\begin{equation}
 \label{eq.ctgammapm}
 \psi _- \coloneq \min _{\theta \in \mathbb{T}} K _0^2(\theta) - \xi
 \ , \qquad \qquad \psi _+ \coloneq \max _{\theta \in \mathbb{T}} K
 _0^2(\theta) + \xi
\end{equation}
with a suitable value of $\xi$.

Let us assume that (depending on $\alpha$ and $\psi _\pm$)
\begin{equation}
\label{eq.boundaryC}
 \partial \mathcal{C} \coloneq  A _\pm \cup B _\pm \cup C _\pm\ ,
\end{equation}
where
\begin{equation*}
 \begin{split}
  A _\pm &\coloneq \{(\theta + \I x, \psi _\pm(\theta) + \I y) \colon
  \theta \in \mathbb{T}, \, |x | \leq \alpha , \, |y| \leq \alpha \}
  \ , \\ B _\pm &\coloneq \{ (\theta \pm \I \alpha ,v + \I w)\colon
  \theta \in \mathbb{T}, \, \psi _-(\theta) \leq v \leq \psi
  _+(\theta), \, |w | \leq \alpha \}\ , \\ C _\pm &\coloneq \{ (\theta
  + \I \sigma ,v \pm \I \alpha)\colon \theta \in \mathbb{T}, \, \psi
  _-(\theta) \leq v \leq \psi _+(\theta), \, |\sigma | \leq \alpha \}
  \ .
 \end{split}
\end{equation*}
We have different cases to get a lower bound on ${\rm dist}(K _0
(\mathbb{T}_{\rho _0}), \partial \mathcal{C})$. Let us consider
generic points
 \begin{align*}
  x &= \bigl(\theta + \I \sigma + \overline{K} _0^1(\theta + \I
  \sigma) \ , \overline{K} _0^2(\theta + \I \sigma)\bigr) \in K _0
  (\mathbb{T}_{\rho _0}) \ , \\ a _1^\pm &= (\theta _1 + \I x _1, \psi
  _\pm(\theta _1) + \I y _1) \in A _\pm \ , \\ b _1^\pm &= (\theta _1
  \pm \I \alpha, v _1 + \I w _1) \in B _\pm \ , \\ c _1^\pm &= (\theta
  _1 + \I \sigma _1, v _1 \pm \I \alpha) \in C _\pm \ .
\end{align*}
If we use $\psi _\pm$ constants, as those defined in
\eqref{eq.ctgammapm}, then we need to compute the $\Upsilon _i$
quantities given by
\begin{equation}
\label{eq.Upsilons}
 \begin{split}
 |x - a _1^+| &\geq \inf _{\theta + \I \sigma \in \mathbb{T}_{\rho _0}} |\Re K _0 ^2 (\theta + \I \sigma) - \psi _+| \eqcolon \Upsilon _1 \ ,\\
 |x - a _1^-| &\geq \inf _{\theta + \I \sigma \in \mathbb{T}_{\rho _0}} |\Re K _0 ^2 (\theta + \I \sigma) - \psi _-| \eqcolon \Upsilon _2 \ ,\\
 |x - b _1^+| &\geq \inf _{\theta + \I \sigma \in \mathbb{T}_{\rho _0}} |\Im K _0 ^1 (\theta + \I \sigma) - \alpha| \eqcolon \Upsilon _3 \ ,\\
 |x - b _1^-| &\geq \inf _{\theta + \I \sigma \in \mathbb{T}_{\rho _0}} |\Im K _0 ^1 (\theta + \I \sigma) + \alpha| \eqcolon \Upsilon _4 \ ,\\
 |x - c _1^+| &\geq \inf _{\theta + \I \sigma \in \mathbb{T}_{\rho _0}} |\Im K _0 ^2 (\theta + \I \sigma) - \alpha| \eqcolon \Upsilon _5 \ ,\\
 |x - c _1^-| &\geq \inf _{\theta + \I \sigma \in \mathbb{T}_{\rho _0}} |\Im K _0 ^2 (\theta + \I \sigma) + \alpha| \eqcolon \Upsilon _6 \ .
 \end{split}
\end{equation}
Thus, if we take $\Upsilon \coloneq \min \{\Upsilon _1, \Upsilon
_2, \Upsilon _3, \Upsilon _4, \Upsilon _5, \Upsilon _6 \}$, then
\begin{equation}
\label{eq.Upsilon}
 {\rm dist}( K _0(\mathbb{T} _{\rho _0}), \partial \mathcal{C}) \geq \Upsilon \ .
\end{equation}
Therefore, we can choose $\zeta$ so that $\Upsilon \geq \zeta > 0$.
Finally, we can set $\Lambda \coloneq (\ecc _0 - \varphi, \ecc _0 +
\varphi)$ with $\varphi \geq \max \{\zeta, 2 \kappa _\ecc \varepsilon
_0\}$ and $\kappa _\ecc$ given in \ref{main.H5}.

Note that the computation of $\Upsilon _i$ in \eqref{eq.Upsilons} does
not need to be rigorous, because we can take $\zeta$ further smaller
than the approximated $\Upsilon$.

A second remark in the computation of $\Upsilon _i$ is that we can use
the complex version of the FFT to make the computation faster. Indeed,
using Appendix~\ref{sec.complexiFou}, we complexify the real
representation of the Fourier coefficients of $K _0$ and then use the
FFT to get the corresponding table of values in an equispaced complex
plane $\mathbb{T} _{\rho _0}$. This process makes the computation of
an approximated $\Upsilon$ efficient, easily running in a today's
laptop without a strong need of concurrency.

\subsection{Steps to approximate the $Q$ quantities} \label{sec:stepsQ}
Once we obtain the initial numerical approximate solution $(K _0, \ecc
_0)$ of the invariance equation \eqref{invariance} via the Newton
Algorithm \ref{alg.newton}, we choose $0 < \rho _0 < 1$ to compute the
different quantities involved in the KAM estimates of the
Theorem~\ref{main}. That means to compute the quantities of
\ref{main.H5} and the constants in Appendix~\ref{app:constants}. The
constants only depend on norms of functions from the
Algorithm~\ref{alg.newton} like $\|D K _0\| _{\rho_0}$, $\| DK
_0^{-1}\| _{\rho_0}$, $\| N \| _{\rho_0}$, $\| S \| _{\rho_0}$,
etc. The $Q$ quantities require more effort and we will use the
procedure described in Section~\ref{sec:Q}.

The first $Q$ quantity $Q _{E _0}$ in Section~\ref{sec:QE0} requires
to choose $0 < \delta _0 < \rho _0$. For the other $Q$ quantities in
Section~\ref{sec:Q} we need first to choose $\xi$ to get $\psi _\pm$
from \eqref{eq.ctgammapm} and $\alpha$ for \eqref{eq.boundaryC}. Then
we compute the $\Upsilon$ such that \eqref{eq.Upsilon} is
satisfied. Finally, we can choose $\zeta$ which is the last crucial
value that fixes all the other quantities to check the inequalities
\eqref{C1}--\eqref{C10}.

We note that the complex quantities of $Q$, which in fact are the
hardest ones, do not need to be extremely rigorous because especially
for those involving high order variational flows, they are always
affected by the multiplication of small values, like the
$\varepsilon_0$, as one can realize looking at \eqref{C1}--\eqref{C10}
and the Appendix~\ref{app:constants}.  Therefore, our approach will
just consider the quantities in a mesh of the six sets in
\eqref{eq.boundaryC}, rather than a rigorous enclosure. In fact, we
also compensate the correctness of our numbers using multiprecision,
that was already needed to reach parameter values close to the
numerical break-down.

\section{KAM estimates for the spin-orbit problem}
\label{sec:estimates}

The application of Theorem~\ref{main} requires to check the conditions
\eqref{C1}--\eqref{C10} that depend on the choice of some parameters.
We did not found a general procedure to select $\rho _0$, $\delta _0$,
$\xi$, $\alpha$ and $\zeta$, so that we can ensure a priori that the
inequalities will be fulfilled. Nevertheless, we provide the values of
these numbers for the cases $\omg _1$ and $\omg _2$ with respective
spin-orbit parameters given in \eqref{eq.conti-example-omg1} and
\eqref{eq.conti-example-omg2}.

In the two cases $\omg _1$ in \eqref{eq.conti-example-omg1} and $\omg
_2$ in \eqref{eq.conti-example-omg2}, by trial and error we have made
the following choice:
\begin{equation}
\label{eq.param-choice}
 \begin{split}
  \rho _0 &= 7.629394531250000 \cdot 10^{-6} = 2^{-17}\ , \\
  \delta _0 &= 9.536743164062500 \cdot 10^{-7} = 2^{-20}\ , \\
  \xi &= 0.0054 \ , \\
  \alpha &= 0.000016\ , \\
  \zeta &= 9.3132257461547851562500 \cdot 10^{-10} = 2^{-30} \ .
 \end{split}
\end{equation}
Note that $\rho _0$, $\delta _0$, and $\zeta$ are just a power of 2,
which means that they have an exact numerical representation in a
computer.

From the choice of values in \eqref{eq.param-choice}, we compute $\psi
_\pm$ and $\Upsilon$ in Section~\ref{sec:boundaryC} using just double
precision:
\begin{equation*}
 \begin{array}{l|cc}
            & \omg _1   & \omg _2 \\ \hline
  \psi_- & \mathtt{2.468595425049463e-01} & \mathtt{2.093861593414215e-01} \\
  \psi_+ & \mathtt{2.682454746721682e-01} & \mathtt{2.306499653402554e-01} \\
  \Upsilon  & \mathtt{1.833012143471895e-06} & \mathtt{1.842114896678543e-06}
 \end{array}
\end{equation*}

\vskip.1in

Then, the $Q$ quantities can be computed following
Section~\ref{sec:stepsQ}. In this computation, we parallelize the
different evaluations in a grid of $16\times 16 \times 16$ points with
a final CPU time of around 33h with 31 threads and 18h with 54
threads. As in the solution computed in Section~\ref{sec.conti}, we
perform all the computations with 170 bits of precision. In
particular, we know that the error in the invariance equation, the
$\varepsilon _0$ in Theorem~\ref{main}, is at most $10^{-45}$ because
it is the requested tolerance in the Newton's process. Moreover, once
all the $Q$ quantities are computed, we perform the final checks of
inequalities \eqref{C1}--\eqref{C10} using a little bit more bits, say
250, to prevent possible overflows in the comparisons.

\vskip.1in

We conclude by saying that the conditions of the theorem are satisfied
for the values given in \equ{eq.conti-example-omg1} for $\omega_1$ and
\equ{eq.conti-example-omg2} for $\omega_2$.  The values of the
quantities needed to prove Theorem~\ref{main} for $\omega_1$ and
$\omega_2$ are listed, respectively, in Appendix~\ref{app:quantities1}
and ~\ref{app:quantities2}.  The values of $\varepsilon$ that we
obtain are essentially coinciding with the numerical break-down
values, which are computed in \cite{CCGL20c}. This result shows the
efficacy of KAM theorem in providing a constructive method to follow
the invariant attractors up to break-down.

\vskip.2in

\appendix
\renewcommand{\thesection}{\S\Alph{section}}
\renewcommand{\thethm}{\Alph{section}.\arabic{thm}}

\section{List of the constants of Theorem~\ref{main}}
\label{app:constants}
The explicit expressions of the constants used in Theorem~\ref{main}
are given below (see \cite{CCL20} for their derivation).

\beqano
C_{\sigma 0}&\coloneq&\mathcal{T}_0\ \Big[|\lambda-1|\ \Big({1\over {||\lambda|-1|}}\|S_0\|_{\rho_0}+1\Big)
+\|S_0\|_{\rho_0}\Big]\ \|M_0^{-1}\|_{\rho_0}\ ,\nonumber\\
C_{W_2 0}&\coloneq&{1\over {||\lambda|-1|}}\Big(1+C_{\sigma 0} Q_{\param }\Big)\|M_0^{-1}\|_{\rho_0}\ ,\nonumber\\
\overline{C}_{W_2 0} &\coloneq& 2\mathcal{T}_0\ \Big({1\over {||\lambda|-1|}}\|S_0\|_{\rho_0}+1\Big)\,
Q_{\param }\ \|M_0^{-1}\|_{\rho_0}^2\ ,\nonumber\\
C_{W_1 0}&\coloneq&C_0\Big(\|S_0\|_{\rho_0}(C_{W_2 0}+\overline{C}_{W_2 0})+\|M_0^{-1}\|_{\rho_0}+Q_{\param} \|M_0^{-1}\|_{\rho_0}C_{\sigma 0}\Big)
\ ,\nonumber\\
C_{W0}&\coloneq&C_{W_1 0}+(C_{W_2 0}+\overline{C}_{W_2 0})\nu\delta_0^{\tau}\ ,\nonumber\\
C_{\eta 0}&\coloneq& C_{W0}\|M_0\|_{\rho_0}+C_{\sigma 0}\nu\delta_0^{\tau}\ ,\nonumber\\
C_{\mathcal{R} 0}&\coloneq&Q_{E 0}(\|M_0\|_{\rho_0}^2 C_{W0}^2+C_{\sigma 0}^2\nu^2\delta_0^{2\tau})\ ,\nonumber\\
C_{{\mathcal E} 0}&\coloneq& C_{W0}\nu \delta_0^{-1+\tau}+C_{\mathcal{R} 0}\ ,\nonumber\\
C_{d0}&\coloneq& C_{W0}\ \|M_0\|_{\rho_0}\ ,\nonumber\\
\kappa_\param&\coloneq&4 C_{\sigma 0}\ ,\nonumber\\
D_K&\coloneq& 4C_{d0}\ \nu^{-1} \delta_0^{-\tau-1}\ \varepsilon_0\ ,\nonumber\\
D_{2K}&\coloneq& 4\ C_{d0} \nu^{-1} \delta_0^{-\tau-2}\ \varepsilon_0\ ,\nonumber\\
C_N&\coloneq&\|N_0\|_{\rho_0}^2\
{{2\|DK_0\|_{\rho_0}+D_K}\over {1-\|N_0\|_{\rho_0}D_K(2\|DK_0\|_{\rho_0}+D_K)}}\ ,\nonumber\\
C_M&\coloneq&1+J_e\Big[C_N(\|DK_0\|_{\rho_0}+D_K)+\|N_0\|_{\rho_0}\Big]\ ,\nonumber\\
C_{Minv}&\coloneq&C_{N}(\|DK_0\|_{\rho_0}+D_K)+\|N_0\|_{\rho_0}+J_e\ ,\nonumber\\
C_S&\coloneq& 2 J_e Q_z \ \Big\{(\|N_0\|_{\rho_0}+C_ND_K)\ \Big[D_K (\|N_0\|_{\rho_0}+C_ND_K)\nonumber\\
&+&\|DK_0\|_{\rho_0}\|N_0\|_{\rho_0}+\|DK_0\|_{\rho_0}C_ND_K\Big]\nonumber\\
&+&C_N\|DK_0\|_{\rho_0} \Big[D_K
(\|N_0\|_{\rho_0}+C_ND_K)
+\|DK_0\|_{\rho_0}\|N_0\|_{\rho_0}+\|DK_0\|_{\rho_0}C_N D_K\Big]\nonumber\\
&+&\|N_0\|_{\rho_0}\|DK_0\|_{\rho_0}(\|N_0\|_{\rho_0}+C_ND_K)
+C_N \|N_0\|_{\rho_0} \|DK_0\|_{\rho_0}^2\Big\}\ ,\nonumber\\
C_{SB}&\coloneq& {1\over
{||\lambda|-1|}}Q_{\param}\|M_0^{-1}\|_{\rho_0}C_S+ 2 J_e Q_z\
\|N_0\|_{\rho_0}^2\ \|DK_0\|_{\rho_0}^2{1\over {||\lambda|-1|}}\
C_{Minv}
\ Q_{\param}\nonumber\\
&+&2C_S\ {1\over {||\lambda|-1|}}\ C_{Minv}\ Q_{\param}\ D_K\ ,\nonumber\\
C_\tau&\coloneq&\max\Big\{C_S, C_{SB}+2C_{Minv} Q_{\param}\Big\}\ D_K\ ,\nonumber\\
C_T&\coloneq& {{\mathcal{T}_0^2}\over {1-\mathcal{T}_0 C_\tau}}\
\max\Big\{C_S, C_{SB}+2C_{Minv} Q_{\param}\Big\}\ ,\nonumber\\
C_\sigma&\coloneq&C_T\ \Big\{|\lambda-1|\ \Big[{1\over {||\lambda|-1|}}(\|S_0\|_{\rho_0}+C_S
D_K)+1\Big]\nonumber\\
&+&\Big(\|S_0\|_{\rho_0}+C_S D_K\Big)\Big\}\ \Big(\|M_0^{-1}\|_{\rho_0}+C_{Minv}D_K\Big)\nonumber\\
&+&\mathcal{T}_0\ \Big\{|\lambda-1|\ \Big[{1\over {||\lambda|-1|}}(\|S_0\|_{\rho_0}+C_S D_K)+1\Big]C_{Minv}\nonumber\\
&+&|\lambda-1|\ {1\over {||\lambda|-1|}}\ \|M_0^{-1}\|_{\rho_0}C_S+C_S
\Big(\|M_0^{-1}\|_{\rho_0}+C_{Minv}D_K\Big)+C_{Minv}\|S_0\|_{\rho_0}\Big\}\ ,\nonumber\\
\overline{C}_{W_2}&\coloneq& 4C_T\ \Big[{1\over
{||\lambda|-1|}}(\|S_0\|_{\rho_0}+C_SD_K)+1\Big]\ Q_{\param}
(\|M_0^{-1}\|_{\rho_0}+D_K)^2\nonumber\\
&+&4\mathcal{T}_0 Q_{\param}\
{1\over {||\lambda|-1|}}C_S\ (\|M_0^{-1}\|_{\rho_0}+D_K)^2\nonumber\\
&+&4\mathcal{T}_0\ Q_{\param} \Big[{1\over {||\lambda|-1|}} (\|S_0\|_{\rho_0}+C_SD_K)+1\Big]
(D_K+2\|M_0^{-1}\|_{\rho_0})\nonumber\\
C_{\mathcal{R}}&\coloneq&Q_{E0}\ \Big[(2C_M\|M_0\|_{\rho_0}+C_M^2D_K)
(C_{W0}+C_WD_K)^2+\|M_0\|_{\rho_0}^2(C_W^2D_K+2C_{W0}\ C_W)\nonumber\\
&+&(C_\sigma^2D_K+2C_{\sigma 0}C_\sigma)\nu^2\delta_0^{2\tau}\Big]
+C_Q\ \Big[(\|M_0\|_{\rho_0}+C_M D_K)^2(C_{W0}+C_W D_K)^2\nonumber\\
&+&(C_{\sigma 0}+C_\sigma D_K)^2\nu^2\delta_0^{2\tau}\Big]\delta_0^{-1}\ ,\nonumber\\
C_{W_2}&\coloneq&{1\over {||\lambda|-1|}}\ \Big[1+2Q_{\param}
\|M_0^{-1}\|_{\rho_0}C_\sigma
+2Q_{\param} C_{\sigma 0}+2Q_{\param}C_\sigma D_K\Big]\ ,\nonumber\\
C_{W_1}&\coloneq& C_0\Big[\|S_0\|_{\rho_0}C_{W_2}+C_S C_{W_2 0}+
C_S C_{W_2}D_K+\|S_0\|_{\rho_0} \overline{C}_{W_2}\nonumber\\
&+&C_S \overline{C}_{W_2 0}+C_S
\overline{C}_{W_2}D_K+1+2Q_{\param}\|M_0^{-1}\|_{\rho_0}C_\sigma+2Q_{\param}
C_{\sigma 0}+
2Q_{\param} C_\sigma D_K\Big]\ ,\nonumber\\
C_W&\coloneq& C_{W_1}+C_{W_2}\nu\delta_0^{\tau} +\overline{C}_{W_2}\nu\delta_0^{\tau}\nonumber\\
C_Q&\coloneq&{1\over 2}\
\max\Big\{1+\sup_{z\in\mathcal{C}}|D^3 f_{\param_0}(z)|\ \|DK_0\|_{\rho_0}^2\delta_0^2\nonumber\\
&+&\sup_{z\in\mathcal{C},\param\in\Lambda,|\param-\param_0|<2\kappa_\param\varepsilon_0}|D_\param D^2 f_{\param}(z)|\ \|DK_0\|_{\rho_0}^2\
{{C_{\sigma 0}}\over {C_{d0}}}\delta_0^{\tau+2}\nonumber\\
&+&\sup_{z\in\mathcal{C}}|D^2 f_{\param_0}(z)|\ \|DK_0\|_{\rho_0}\ \delta_0\nonumber\\
&+&\sup_{z\in\mathcal{C}}|D^3 f_{\param_0}(z)|\ \|DK_0\|_{\rho_0}\ 4 C_{d0} \nu^{-1}\delta_0^{-\tau+1}\varepsilon_0\nonumber\\
&+&\sup_{z\in\mathcal{C},\param\in\Lambda,|\param-\param_0|<2\kappa_\param\varepsilon_0}|D_\param D^2 f_{\param}(z)|\ \|DK_0\|_{\rho_0}\
4 C_{\sigma 0} \delta_0\varepsilon_0\nonumber\\
&+&\sup_{z\in\mathcal{C}}|D^2 f_{\param_0}(z)|\ \|D^2 K_0\|_{\rho_0}^2\ \delta_0^2\nonumber\\
&+&\sup_{z\in\mathcal{C},\param\in\Lambda,|\param-\param_0|<2\kappa_\param\varepsilon_0}|D_\param D f_{\param}(z)|\ \|D^2K_0\|_{\rho_0}\
{{C_{\sigma 0}}\over {C_{d0}}}\nu\delta_0^{\tau+2}\nonumber\\
&+&\sup_{z\in\mathcal{C}}|D^2 f_{\param_0}(z)|\ (\|DK_0\|_{\rho_0}+D_K)\ \delta_0\nonumber\\
&+&\sup_{z\in\mathcal{C}}|D^3 f_{\param_0}(z)|\ (\|DK_0\|_{\rho_0}+D_K)4 C_{d0}\ \nu^{-1}\delta_0^{-\tau+1}\varepsilon_0\nonumber\\
&+&\sup_{z\in\mathcal{C},\param\in\Lambda,|\param-\param_0|<2\kappa_\param\varepsilon_0}|D_\param D^2 f_{\param}(z)|\ (\|DK_0\|_{\rho_0}+D_K)\
4C_{\sigma 0}  \delta_0\varepsilon_0\nonumber\\
&+&\sup_{z\in\mathcal{C}}|D f_{\param_0}(z)|+\sup_{z\in\mathcal{C}}|D^2 f_{\param_0}(z)|\ (4 C_{d0}\nu^{-1}\delta_0^{-\tau})\ \varepsilon_0\nonumber\\
&+&\sup_{z\in\mathcal{C},\param\in\Lambda,|\param-\param_0|<2\kappa_\param\varepsilon_0}|D_\param D f_{\param}(z)|\ \kappa_\param\varepsilon_0\ ,\nonumber\\
&&\sup_{z\in\mathcal{C}}|DD_\param f_{\param_0}(z)|\ \delta_0+
\sup_{z\in\mathcal{C}}|D^2 D_\param f_{\param_0}(z)|\ \delta_0^2\ (\|DK_0\|_{\rho_0}+D_K)\nonumber\\
&+&\sup_{z\in\mathcal{C},\param\in\Lambda,|\param-\param_0|<2\kappa_\param\varepsilon_0}|DD_\param^2 f_{\param}(z)|\
{{C_{\sigma 0}}\over {C_{d0}}}\nu\delta_0^{\tau+2}\ (\|DK_0\|_{\rho_0}+D_K),\nonumber\\
&&\sup_{z\in\mathcal{C},\param\in\Lambda,|\param-\param_0|<2\kappa_\param\varepsilon_0}|D_\param^3 f_{\param}(z)|\
{{C_{\sigma 0}}\over {C_{d0}}}\nu\delta_0^{\tau+2}\Big\}\ .\nonumber\\
\eeqano

\section{Newton's algorithm}
\label{app:algorithm}

In this Section, we provide Newton's algorithm for finding an
invariant attractor of the spin-orbit problem; the algorithm is fully
detailed in \cite{CCGL20a}.

\begin{alg}[Newton's method for finding a torus in the spin-orbit problem]
 \label{alg.newton}
 \
 \begin{enumerate}
\setlength{\itemsep}{.7em}
\renewcommand*{\theenumi}{\emph{\arabic{enumi}}}
\renewcommand*{\labelenumi}{\theenumi.}
  \item [$\star$] \texttt{Inputs:} A fixed frequency $\omega$, the
    conformally symplectic map $P _\ecc$ given in \eqref{eq.Pe} for
    fixed values of the parameters $\eps$ and $\dis$. Initial values of the
    unknowns; the eccentricity $\ecc$ and the embedding $K\colon
    \mathbb{T} \rightarrow \mathbb{T} \times \mathbb{R}$.
  \item [$\star$] \texttt{Output:} New $K$ and $\ecc$ satisfying the
    invariance equation \eqref{invariance} up to a given tolerance.
  \item [$\star$] \texttt{Notation:} If $A$ is a function defined in
    $\mathbb{T}$, $\overline A \coloneq \int _{\mathbb{T}} A$ and $A^0
    \coloneq A - \overline A$.
  \item \label{alg.newton-error} $E \gets P _ \ecc \circ K - K \circ T
    _\omega$ denote the components $E \coloneq (E _1, E _2)$, \newline
    $E _1 \gets E _1 - \mathtt{round}(E _1)$.
  \item $\alpha \gets D K$.
  \item $N \gets (\alpha ^t \alpha) ^{-1}$.
  \item \label{alg.newton-M} $M \gets
  \begin{bmatrix}
   \alpha & J ^{-1} \alpha N
  \end{bmatrix}$.
  \item $\widetilde E \gets (M ^{-1}\circ T _\omega) E$.
  \item $\lambda$ given in \eqref{lambda}.
  \item \label{alg.newton-PSA} $P \gets \alpha N$, \\ $S \gets (P
    \circ T _\omega)^t D P _ \ecc \circ K J ^{-1} P$, \\ $\widetilde A
    \gets M ^{-1} \circ T _\omega D _\ecc P _\ecc \circ K$ denote the
    components $\widetilde A \coloneq (\widetilde A _1, \widetilde A
    _2)$.
  \item \label{alg.newton-cohom1} $(B _a)^0$ solving $\lambda (B _a)
    ^0 - (B _a)^0 \circ T _ \omega = - (\widetilde E _2)^0$, \\ $(B
    _b)^0$ solving $\lambda (B _b) ^0 - (B _b)^0 \circ T _ \omega = -
    (\widetilde A _2)^0$.
  \item Find $\overline W _2$, $\sigma$ solving the linear system
  \begin{equation*}
   \begin{pmatrix}
    \overline {S} & \overline{S(B _b)^0} +
\overline{\widetilde A _1} \\
    \lambda - 1 &  \overline{\widetilde A _2}
   \end{pmatrix}
   \begin{pmatrix}
     \overline W _2 \\ \sigma
   \end{pmatrix} =
   \begin{pmatrix}
     - \overline{\widetilde E _1} - \overline{S(B _a)^0} \\
     - \overline{\widetilde E _2}
   \end{pmatrix}
  \end{equation*}
  \item $(W _2)^0 \gets (B _a)^0 + \sigma (B _b)^0$.
  \item $W _2 \gets (W _2)^0 + \overline{W}_2$.
  \item \label{alg.newton-cohom2}
  $(W _1)^0$ solving $(W _1) ^0 - (W _1)^0 \circ T _ \omega = - (S W _2)^0
- (\widetilde E _1)^0 - (\widetilde A _1)^0 \sigma$.
  \item $K \gets K + M W$, \\
  $\ecc \gets \ecc + \sigma$.
  \item Iterate from \eqref{alg.newton-error} until convergence in $E$
    with a prescribed tolerance $\tilde\epsilon$.
 \end{enumerate}
\end{alg}

\section{Complexification of a Fourier series}
\label{sec.complexiFou}

If $x \colon \mathbb{T} \rightarrow \mathbb{R}$ is a periodic and
smooth mapping of period 1, it admits an $N$-th order truncated
Fourier series with Fourier coefficients $\{ x _k \} _{k = 0}^{N-1}
\subset \mathbb{R}$:
\begin{equation}
\label{eq.Foucoefs}
 x(\theta) = \frac{x _0}{2} + \frac{x _{N/2}}{2} \cos (\pi N \theta) +
 \sum_{k = 1}^{N/2-1} x _{2k} \cos (2\pi k \theta) + x _{2k+1}
 \sin(2\pi k \theta)\ .
\end{equation}
For simplicity and easy notation we assume $N$ to be an even positive
integer in \eqref{eq.Foucoefs}. The complexification process of the
map $x$ consists in lifting the spaces $\mathbb{T}$ and $\mathbb{R}$
to the complex numbers such that it coincides with $x$ when it is
restricted to the real values.

To make it simpler, it is convenient to extend the quantity of real
numbers in \eqref{eq.Foucoefs} and make explicit the symmetry in the
complex version. In other words, \eqref{eq.Foucoefs} is equivalent to
\begin{equation*}
 x (\theta) = x _0 + 2\sum_{k = 1}^{N/2} (x _{2k} - \I x _{2k+1}) e^{2 \pi k \I \theta} + (x _{2k} + \I x _{2k+1}) e^{-2 \pi k \I \theta}
\end{equation*}
with $x _{N+1} = 0$. Now, if $\rho > 0$, then
\begin{equation*}
 x(\theta + \I \rho) = x _0 + 2\sum_{k = 1}^{N/2} (x _{2k} - \I x _{2k+1}) e^{2 \pi k \I (\theta + \I \rho)} + (x _{2k} + \I x _{2k+1}) e^{-2 \pi k \I (\theta + \I \rho)},
\end{equation*}
which allows one to provide the Fourier coefficients $\{ (x _{2k}
\pm \I x _{2k+1}) e^{\pm 2\pi k \rho}\}$ making the initial real
Fourier expression to a complex one.

\section{KAM quantities for the frequency $\omg _1$}
We list below the quantities needed to implement Theorem~\ref{main} to
get the existence of an invariant attractor with frequency $\omega_1$.

\label{app:quantities1}
 \renewcommand{\arraystretch}{1.3}
 \begin{longtable}[l]{rcl}
  $N$ &=& $16384$ \ , \\
$\eps$ &=& $\mathtt{1.1632963641877116367716112642948530559675531382297e-02}$ \ , \\
$\dis$ &=& $10^{-3}$ \ , \\
$\ecc$ &=& $\mathtt{3.1675286891174832107186084513865661761571784973618e-01}$ \ , \\
$\rho _0$ &=& $2^{-17}$ \ , \\
$\delta _0$ &=& $2^{-20}$ \ , \\
$\|DK\|_{\rho_0}$ &=& $\mathtt{6.2076969839032564048438325650777912419214845002300e+00}$ \ , \\
$\|DK^{-1}\|_{\rho_0}$ &=& $\mathtt{1.8328129957258449874460075408233923038434712690096e+05}$ \ , \\
$\|D^2K\|_{\rho_0}$ &=& $\mathtt{1.1686089945113448858821651745887573374665126081719e+02}$ \ , \\
$Q _{E _0}$ &=& $\mathtt{1.9132315264792576102165122788680383078808432879626e+00}$ \ , \\
$\|N\|_{\rho_0}$ &=& $\mathtt{9.8051171808495981670035137469708799108365949325248e+00}$ \ , \\
$\|N^{-1}\|_{\rho_0}$ &=& $\mathtt{1.0113946410899826227827056006594783412357401114959e+01}$ \ , \\
$\|S\|_{\rho_0}$ &=& $\mathtt{5.7223321830936249091412643788103653938262105245420e+01}$ \ , \\
$\|E _0\|_{\rho_0}$ &=& $\mathtt{5.7356559781857403764979281930553140398186337716656e-48}$ \ , \\
$\lambda$ &=& $\mathtt{9.8689359923042965027116069623508749107899367134535e-01}$ \ , \\
$\|M\|_{\rho_0}$ &=& $\mathtt{1.2040958250027560141817737598227413458417075907998e+01}$ \ , \\
$\|M^{-1}\|_{\rho_0}$ &=& $\|M\|_{\rho_0}$ \ , \\
${\mathcal T}_0$ &=& $\mathtt{9.9819949009440259228900924748534932771289016437641e+01}$ \ , \\
$8C_\sigma\|E _0\|_{\rho_0}$ &=& $\mathtt{6.3125418322117269519458608993574236031175483317508e-42}$ \ , \\
$\zeta$ &=& $2^{-30}$ \ , \\
$Q_z$ &=& $\mathtt{6.6101300016209423195423975547239258176432851802452e+00}$ \ , \\
$Q_{\param}$ &=& $\mathtt{1.4175899711779293156363275537004756604799008098606e-01}$ \ , \\
$Q_{zz}$ &=& $\mathtt{2.7720843711101391648970926156205952547296169675164e+01}$ \ , \\
$Q_{\param z}$ &=& $\mathtt{1.8953747809385677544688739634954194030662439685438e+00}$ \ , \\
$Q_{zzz}$ &=& $\mathtt{1.0724765398115262597036815561828324853342846787822e+03}$ \ , \\
$Q_{\param zz}$ &=& $\mathtt{3.5749217993717541990122718539427749511142822626074e+02}$ \ , \\
$Q_{z\param }$ &=& $\mathtt{1.8953747809385677544688739634954194030663082259039e+00}$ \ , \\
$Q_{\param\param}$ &=& $\mathtt{5.0662817916743259572206058032898398336985496305233e-01}$ \ , \\
$Q_{zz\param}$ &=& $\mathtt{3.5749217993717541990122718539427749511143316805187e+02}$ \ , \\
$Q_{\param\param z}$ &=& $\mathtt{6.5647293504204776520816015502477982907961893178051e+01}$ \ , \\
$Q_{\param\param\param}$ &=&
$\mathtt{1.1476657592536890159871266106814665145940757444159e+00}$
\ .
 \end{longtable}

%

\section{KAM quantities for the frequency $\omg _2$}
We list below the quantities needed to implement
Theorem~\ref{main} to get the existence of an invariant attractor
with frequency $\omega_2$.

 \label{app:quantities2}
 \begin{longtable}[l]{rcl}
  $N$ &=& $4096$ \ , \\
$\eps$ &=& $\mathtt{1.2697630024415883032123830013667613509009950826168e-02}$ \ , \\
$\dis$ &=& $10^{-3}$ \ , \\
$\ecc$ &=& $\mathtt{2.4824740823563165902227100091869770425731996450084e-01}$ \ , \\
$\rho _0$ &=& $2^{-17}$ \ , \\
$\delta _0$ &=& $2^{-20}$ \ , \\
$\|DK\|_{\rho_0}$ &=& $\mathtt{6.2401368092989368560939911390480948796213323884872e+00}$ \ , \\
$\|DK^{-1}\|_{\rho_0}$ &=& $\mathtt{9.7663343052106062599854114524341354648300957997991e+04}$ \ , \\
$\|D^2K\|_{\rho_0}$ &=& $\mathtt{1.2599262190633202679574003877751236478676849823924e+02}$ \ , \\
$Q _{E _0}$ &=& $\mathtt{3.7283183855924988259949473978598408908275342300333e+00}$ \ , \\
$\|N\|_{\rho_0}$ &=& $\mathtt{9.7219870102188805011710709653101075119387149314734e+00}$ \ , \\
$\|N^{-1}\|_{\rho_0}$ &=& $\mathtt{1.0224486155736666017813494196253391421224297676866e+01}$ \ , \\
$\|S\|_{\rho_0}$ &=& $\mathtt{5.6566290718009094885071045850417592994899924965064e+01}$ \ , \\
$\|E _0\|_{\rho_0}$ &=& $\mathtt{4.5110963829895625372478056855241916107240582063354e-45}$ \ , \\
$\lambda$ &=& $\mathtt{9.9012510148807761346816298772561891586174978261238e-01}$ \ , \\
$\|M\|_{\rho_0}$ &=& $\mathtt{1.2040013601997889491301308242283364695245420720597e+01}$ \ , \\
$\|M^{-1}\|_{\rho_0}$ &=& $\|M\|_{\rho_0}$ \ , \\
${\mathcal T}_0$ &=& $\mathtt{6.0557474279802520066531787357919583737990862560932e+01}$ \ , \\
$8C_\sigma\|E _0\|_{\rho_0}$ &=& $\mathtt{2.9770931760274406778788288754482772991065644242739e-39}$ \ , \\
$\zeta$ &=& $2^{-30}$ \ , \\
$Q_z$ &=& $\mathtt{6.5592165251990406445369341061622617126571276578225e+00}$ \ , \\
$Q_{\param}$ &=& $\mathtt{1.5083817512231203986293089732663386692582367665897e-01}$ \ , \\
$Q_{zz}$ &=& $\mathtt{2.7092396727127668081914126144670011929964472351248e+01}$ \ , \\
$Q_{\param z}$ &=& $\mathtt{2.7606921497436169824355915345916507538567998407721e+00}$ \ , \\
$Q_{zzz}$ &=& $\mathtt{1.0002777586041620153665189720104193672993532271333e+03}$ \ , \\
$Q_{\param zz}$ &=& $\mathtt{3.3342591953472067178883965733680645576645107571111e+02}$ \ , \\
$Q_{z\param}$ &=& $\mathtt{2.7606921497436169824355915345916507538735860805555e+00}$ \ , \\
$Q_{\param\param}$ &=& $\mathtt{2.8395238802380805094507234385691115589067408070788e-01}$ \ , \\
$Q_{zz\param}$ &=& $\mathtt{3.3342591953472067178883965733680645576942037378572e+02}$ \ , \\
$Q_{\param\param z}$ &=& $\mathtt{7.2062924226872994236028388252435296803169026567968e+01}$ \ , \\
$Q_{\param\param\param}$ &=&
$\mathtt{5.7795906592094838240953693819846769451926848120880e-01}$
\ .
 \end{longtable}

\section{Multivariate polynomials of degree 2}
\label{sec.multipoly2}
Let us consider a polynomial with $d$ variables and degree $2$, namely
\begin{equation}
\label{eq.polydeg2}
 p(x) = p _0 + \sum _{|k| = 1} p _k x^k + \sum _{|k|=2} p _k x^k,
 \qquad k \in \mathbb{N}^d, x = (x _0, \dotsc, x _{d-1})
\end{equation}
with the multi-index conventions $|k| = k _0 + \dotsb + k _{d-1}$ and
\begin{equation}
 x ^k = x _0 ^{k _0} x _1^{k _1} \dotsb x _{d-1} ^{k _{d-1}}.
\end{equation}
Note that in the case of degree $2$, the multi-index $k$ can be
encoded with the canonical vector ${\underline e} _l = (0, \dotsc,
1, \dotsc, 0)$ with $0 \leq l < d$ in $\mathbb{R}^d$. That is,
either ${\underline e} _i$ for $|k|=1$ or ${\underline e} _i +
{\underline e} _j$ with $i \geq j$ for $|k|=2$.

Let us now define $\chi (i) = \#\{ k \in \mathbb{N}^i \colon |k|=2\}$,
which is computable by the recurrence
\begin{equation}
\label{eq.chi}
 \begin{split}
  \chi(0)&=0, \\
  \chi(i) &= \chi(i-1) + i, \quad i \geq 1.
 \end{split}
\end{equation}
Thus the number of elements to store in a computer for
\eqref{eq.polydeg2} is $\chi(d)+d+1$.

The crucial operation for an arithmetic of elements like
\eqref{eq.polydeg2} is the product, in which the key step is the
product of the two homogenous polynomials of degree 1, since the other
terms are just multiplications by the independent term of each of the
polynomials involved. To this end, we must fix a monomial order to
encode the physical index of each of the monomials of degree
$2$. Among all of them, we consider the reverse lexicographical order,
which is illustrated in Table~\ref{tab.monomials} up to $5$ variables.

\begin{table}[ht]
 \[
  \begin{array}{c|r|c|c}
   &\text{Monomial} & \text{Multi-index}  & \text{Index} \\ \hline \hline
   \chi(1)=1
   & x _0^2    & \ve e _0 + \ve e _0 & 0   \\ \hline
   \multirow{2}{*}{$\chi(2)=3$}
   & x _1 x _0 & \ve e _1 + \ve e _0 & 1   \\
   & x _1^2    & \ve e _1 + \ve e _1 & 2   \\ \hline
   \multirow{3}{*}{$\chi(3)=6$}
   & x _2x _0  & \ve e _2 + \ve e _0 & 3   \\
   & x _2 x _1 & \ve e _2 + \ve e _1 & 4   \\
   & x _2^2    & \ve e _2 + \ve e _2 & 5   \\ \hline
   \multirow{4}{*}{$\chi(4)=10$}
   & x _3x _0  & \ve e _3 + \ve e _0 & 6   \\
   & x _3 x _1 & \ve e _3 + \ve e _1 & 7   \\
   & x _3 x _2 & \ve e _3 + \ve e _2 & 8   \\
   & x _3^2    & \ve e _3 + \ve e _3 & 9   \\ \hline
   \multirow{5}{*}{$\chi(5)=15$}
   & x _4x _0  & \ve e _4 + \ve e _0 & 10  \\
   & x _4 x _1 & \ve e _4 + \ve e _1 & 11  \\
   & x _4 x _2 & \ve e _4 + \ve e _2 & 12  \\
   & x _4 x _3 & \ve e _4 + \ve e _3 & 13  \\
   & x _4^2    & \ve e _4 + \ve e _4 & 14
  \end{array}
 \]
 \caption{Bijection encoding between the exponent $x^k$ with $x=(x _0,
   \dotsc, x _{d-1})$, $k = \ve e _i + \ve e _j$, $i \geq j$, where
   $\ve e _l = (0,\dotsc, 1, \dotsc 0)$ and the location on the array
   containing the terms $p _k$ of \eqref{eq.polydeg2}.}
 \label{tab.monomials}
\end{table}

Thus, the location in the array corresponding to $k = {\underline e}
_i + {\underline e} _j$ with $i \geq j$ is given by $\chi(i) + j$. We
implement this procedure in the function \verb"ex2pl(i,j)" given
below.

On the other hand, to know the $i$ and $j$ for a given index $l$ in
the vector of coefficients, one first performs a binary search to know
$k$ such that $\chi(k) \leq l < \chi(k+1)$, then $i = k$ and $j =
l-\chi(k)$.

A possible pseudo code to compute the product $pq$ of two homogeneous
polynomials $p$ and $q$ with $d$ variables and of degree $1$ can then
be
\begin{verbatim}
    int ex2pl(i,j): return chi(max(i,j)) + min(i,j)

    void php1(d,p,q,flag,pq):
        if (flag==0) for (i = 0; i < chi[d]; i++) pq[i]=0
        for (i = 0; i < d; i++) for (j = 0; j < d; j++)
          pq[ex2pl(i,j)]+= p[i] * q[j]
\end{verbatim}

\vskip.1in

Once the product of multivariate polynomials of degree $2$ is clear,
the other elementary operations such as division, power, trigonometric
operations and hyperbolic trigonometric operations can be derived in a
recurrence manner, see \cite{Haro2016}. For instance, the division of
$r(x) = p(x)/q(x)$ has the following terms
\begin{equation}
 \begin{aligned}
  r _0 &= \frac{p _0}{q _0}, \\ r _k &= \frac{p _k - r _0 q _k }{q _0}
  , & |k|&=1\ , \\ r _k &= \frac{1}{q _0} \biggl[p _k - r _0 q _k -
    \underbrace{\biggl(\sum _{|j|=1} q _j x^j \biggr)\biggl(\sum
      _{|j|=1} r _j x^j \biggr)} _{\text{call to the \texttt{php1}
        function}}\biggr], & |k|&=2\ .
 \end{aligned}
\end{equation}


\newcommand{\etalchar}[1]{$^{#1}$}
\def\cprime{$'$} \def\cprime{$'$} \def\cprime{$'$} \def\cprime{$'$}

\end{document}